\documentclass{article}
\pdfoutput=1

\usepackage[utf8]{inputenc}
\usepackage[margin = 1in]{geometry}
\usepackage{amsmath}
\usepackage{amsfonts}
\usepackage{amssymb}
\usepackage{overpic}
\usepackage{epstopdf}
\usepackage{graphicx}
\usepackage{bm}
\usepackage{caption}
\usepackage{subcaption}
\usepackage{multirow}
\usepackage{authblk}
\usepackage[UKenglish]{datetime}
\newdateformat{UKvardate}{%
	\THEDAY\ \monthname[\THEMONTH] \THEYEAR}

\newcommand{\dd}{\ensuremath{\mathrm{d}}}
\newcommand{\Ord}[1]{\ensuremath{\mathcal{O}\left( #1\right)}}
\newcommand{\ns}[1]{\ensuremath{\uppercase{ #1 }}}
\newcommand{\sign}[1]{\ensuremath{\mathrm{sign}\left( #1 \right)}}
\renewcommand{\vec}[1]{\ensuremath{\bm{#1}}}	
\newcommand{\p}{\ensuremath{\bm{p}}}
\newcommand{\s}{\ensuremath{\bm{s}}}
\newcommand{\z}{\ensuremath{\bm{z}}}
\newcommand{\Da}{\ensuremath{\mathcal{D}_{att.}}}
\newcommand{\Dr}{\ensuremath{\mathcal{D}_{rep.}}}
\newcommand{\Dc}[1]{\ensuremath{\mathcal{D}_{C^{ #1 }}}}

\newcommand{\ve}{\ensuremath{\varepsilon}}
\newcommand{\eps}{\ensuremath{\varepsilon}}
\newcommand{\ga}{\ensuremath{\alpha}}
\newcommand{\gk}{\ensuremath{\kappa}}
\newcommand{\gd}{\ensuremath{\delta}}
\newcommand{\gu}{\ensuremath{\upsilon}}
\newcommand{\gz}{\ensuremath{\zeta}}
\newcommand{\gt}{\ensuremath{\theta}}
\newcommand{\gp}{\ensuremath{\varphi}}
\newcommand{\gv}{\ensuremath{\nu}}
\newcommand{\nsomega}{\ensuremath{\Omega}}

\title{Analysis of point-contact models of the bounce of a hard spinning ball 
on a compliant frictional surface.}
\author[$\dagger$,1]{Stanis{\l}aw W.~Biber}
\author[1]{Alan R.~Champneys}
\author[1]{Robert Szalai}
\affil[1]{Department of Engineering Mathematics, University of Bristol, BS8 
1TW, UK} 
\affil[$\dagger$]{Corresponding author: s.biber@bristol.ac.uk}
\UKvardate
\date{\today}

\begin{document}

\maketitle

\begin{abstract}
Inspired by the turf-ball interaction in golf, this paper seeks to 
understand 	the bounce of a ball that can be modelled as a rigid sphere and the 
surface as supplying an elasto-plastic contact force in addition to Coulomb 
friction. A general formulation is proposed that models the finite time 
interval of bounce from touch-down to lift-off. Key to the analysis is 
understanding transitions between slip and roll during the bounce. Starting 
from the rigid-body limit with a an energetic or Poisson coefficient of 
restitution, it is shown that slip reversal during the contact phase 
cannot be captured in this case, which result generalises to the case of 
pure normal compliance. Yet, the introduction of linear tangential 
stiffness and damping,  does enable slip reversal. 
This result is extended to general weakly nonlinear normal and tangential 
compliance. An analysis using Filippov theory of piecewise-smooth systems 
leads to an argument in a natural limit that lift-off while rolling is 
non-generic and	that almost all trajectories that lift off, do so under slip 
conditions. Moreover, there is a codimension-one surface in the space of 
incoming velocity and spin which divides balls that lift off with backspin from 
those that lift off with topspin.  The results are compared with recent 
experimental measurements on golf ball bounce and the theory is shown to 
capture the main features of the data.

\end{abstract}

\begin{center}
\end{center}

\maketitle

\section{Introduction} \label{sec:intro}
Golf is a highly technical sport with large amounts of sponsorship and
prize money. Yet, until relatively recently there have been relatively
few academic studies on the dynamics of a golf ball. 

Those studies which do exist tend to focus on club-ball interaction and on the  
aerodynamic properties of the ball in flight. Perhaps the most influential work 
is the paper by Quintavalla \cite{quintavalla_flight}, which provides a 
practical, 
parametrised model for the flight of a golf ball.
Backed by a physical principles, the model relies on relatively few input
parameters, such as lift and drag coefficients.
The ease with which these can be estimated from experimental campaigns, have 
made the model one of the most well-known and commonly used in the field, both 
by manufacturers and legislators of the game of golf, and has enabled this 
mathematical model to define industry standards. 

Following its flight however, a golf ball will bounce, perhaps several times, 
then typically roll, before coming to rest. Surprisingly 
little has been established for this bounce and roll phase. One 
of the difficulties is that while the launch conditions of a golf shot are 
typically well controlled, those of its landing are not. There are many 
variables both related to the surface and viscoelastic properties of the turf, 
as well the spin and velocity of the ball upon first landing. 

Current literature on ball bounce tends to focus on 
cases where a deformable ball bounces on a rigid surface -- see for example
\cite{Daish,brake_analytical,cricket_impact,cordingley_phd,ghaednia_review}.
The results of such studies are more applicable to sports such as 
tennis, cricket or football, but less is known for the case of an almost
rigid ball  bouncing on a compliant surface, the type of bounce most commonly
observed in golf.

One of the most comprehensive studies of golf ball-turf interaction
was that of Haake \cite{haake_apparatus}, with the aim of
identifying apparatus for a quantitative classification of golf turf.
A number of models of bounce were also
considered as 
part of that investigation, involving nonsmooth transitions between contact and 
non-contact phases, and also between slip and stiction.
However, the models contained a wide range of parameters that would be hard
to measure in the field, thus limiting their
predictive power. 
A limited experimental data 
set was available at the time and as such the model never became widely 
verified and used.

Based on the same data set as Haake, a new study by Penner \cite{penner2002} 
introduced an idea where the multiple forces acting on the ball's surface 
during a compliant bounce are summed to a single force acting at a point at an 
angle. Penner thus suggested that a compliant bounce can be modelled as a rigid 
bounce (as defined in \cite{Daish}) against an inclined surface. This appeared 
to agree with the limited data set available, but appears to be too simplistic 
for a wider range of initial conditions.

Another study looking at golf ball bounce was presented by Cross in 
\cite{cross} and models the bounce using a generalisation of the concept of
Newtonian rigid impact to include a tangential coefficient of restitution
(see also \cite{Roh}). Such a phenomenological
coefficient can be fit to data and is said to allow for
backwards bounce, but there is little physics incorporated into such a 
parameter and it cannot be used to predict {\it a priori} what happens at the 
transition between
roll and slip occur during a bounce.  Once again, the theoretical result is 
backed by  a limited set of experimental data. 

Although each of the studies considered and explained a wide range of 
behaviours seen in golf (such as spin reversal) the lack of 
experimental validation or sound physical background leaves a gap in the 
applicability of the models. We thus seek a model that can be physically 
valid,  can be described using a minimal 
number of parameters, and can be verified experimentally. 

The key issue seems to be to understand possible transitions between slip and 
roll during the  bounce. Such an approach requires the definition of a dry 
friction model; which can be a painful process to match to data -- a central 
problem in the field of tribology \cite{Barber}.
Common dry friction models include the simple Coulomb model, or its 
generalisation to models such as those due to 
Stribeck \cite{stribeck}, LuGre \cite{lugre}, Dahl \cite{dahl}. The state of 
the art seems to 
be so-called 
rate-and-state friction which includes effects of slow creep and pre-slip 
through an auxiliary state variable that tries
to capture the deformation of surface asperities; see e.g.~\cite{Putelat} and 
references therein. Fitting such models' parameters to data can be crucial 
during sustained contact motion.
However, in this work, we shall stick to the simple Coulomb model, as what 
seems to be important is understanding the nature of the transition between 
slip and stick (roll) during a rapid bounce. Also, we shall allow for a change 
in the nature of the surface during the bounce through the tangential 
compliance of the surface, rather than through more complex friction laws. We 
shall also allow for quite general forms of coupling between tangential and
normal compliance, as in \cite{Hoffmann} under certain natural scaling 
hypotheses on the relative size of normal and tangential forces. 

\begin{figure}
	\centering
	\includegraphics[width = 0.6\linewidth]{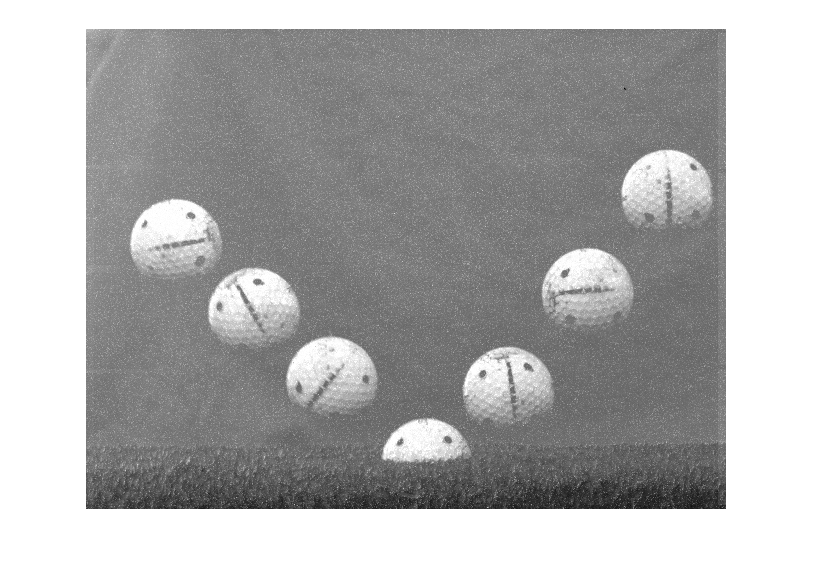}
	\caption{A combination of frames from a recorded high-speed video of a ball 
		bounce. }
	\label{fig:vid_still}
\end{figure}

Modelling processes with rigid bounce and Coulomb friction can be a complex 
process
and is known to lead to several different inconsistencies that are loosely 
termed the
Painlev\'{e} paradox \cite{Varkonyi,Hogan,varkonyi2}. However such paradoxes 
only occur for
slender objects (where there is a sufficiently large coupling between normal 
and tangential degrees of freedom at contact) and do not come into play in this 
work. As we shall see, though, there are a range of finite-time singularities 
that are inherent in nonsmooth frictional contact problems associated with 
transitions between stick and slip, as are inherent in this work.
These can be understood in 
terms of the theory of Filippov systems, see 
\cite{filippov,piecewisebook,jeffrey_hidden} and references therein. 

In this paper we analyse the behaviour of a spinning uniform sphere as it 
bounces off a 
generic nonlinear compliant surface. Although applicable to a wide range of 
problems, our physical constraints and intuition will be directed here by the 
example of a golf ball bouncing. The key question we want to address is 
whether, under parsimonious assumptions about the nature of the surface 
presented in detail in Sec.~\ref{sec:analysis}, equations 
\eqref{eq:functions_constraints}-\eqref{eq:functions_constraints_horizontal_scaling},
a ball that enters a bounce spinning and subsequently grips the surface and 
enters a roll phase during the bounce, will lift off rolling or will lift off 
with reverse spin.

The rest of the paper is outlined as follows. We begin in Section
\ref{sec:preliminaries} with a summary of the problem and
present some experimental observations the details
of which will be presented elsewhere.  We show how a purely
rigid-bounce model cannot account for what is observed, in which
the ball is found to lift off with either topspin or backspin, but
rarely, if ever, in a state of rolling.  An improved model is
presented in Section \ref{sec:two_cases} which captures the finite
time of contact and allows for linear visco-elastic behaviour in both
the normal and tangential degrees of freedom. It is shown that such a
model gives the correct qualitative behaviour.  Section
\ref{sec:analysis} then presents a generalised visco-elastic model,
that allows for nonlinearity and coupling between normal and
tangential degrees of freedom. We analyse the critical transition
between topspin and backspin at lift-off and show that this
represents a so-called two-fold singularity within the Filippov
systems.  We conclude that, in the natural limit of a small ratio
between tangential and normal stiffness, the lift-off with rolling
represents the codimension-one manifold of trajectories that passes
through the singularity. Moreover, this manifold separates two open
sets of initial conditions that enter a rolling state at some stage
during the bounce, yet lift off with topspin, or with backspin.
Section \ref{sec:conclusion} contains conclusions and discusses open
problems.

\section{Preliminaries}
\label{sec:preliminaries}

\begin{figure}[!htb]
	\centering
	\centering
	\begin{overpic}[width = 0.4\linewidth]{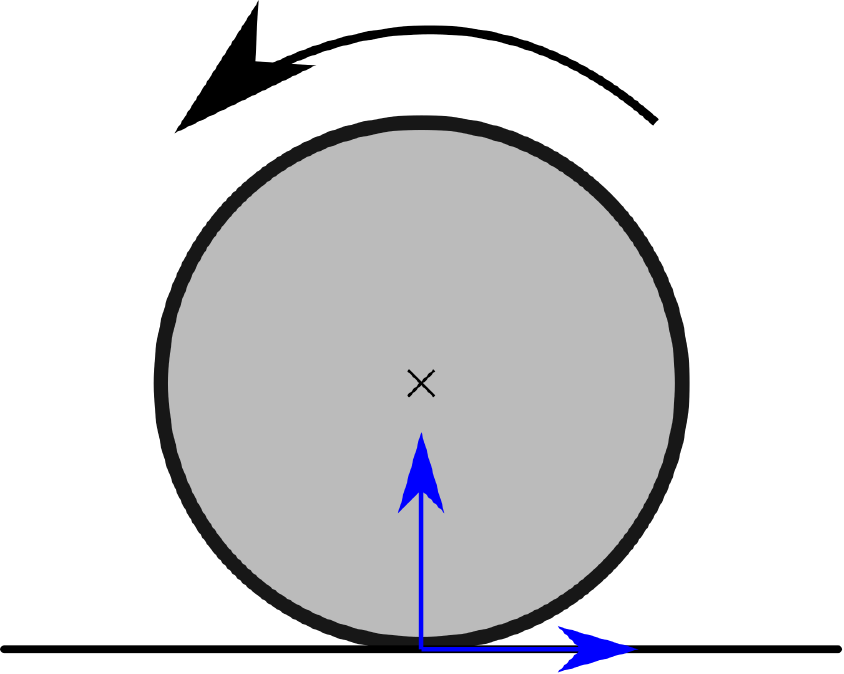}
		\put(75,5){ $\lambda_T$}
		\put(55,20){ $\lambda_N$}
		\put(45,40){ $(x,y)$}
		\put(78,80){ $\omega$}
	\end{overpic}
	\caption{Setting of the notation used. }
	\label{fig:setting_notation}
\end{figure}
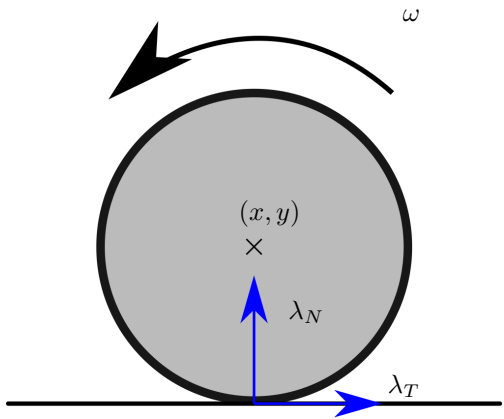

Throughout this paper, we consider an isotopic rigid sphere of uniform density 
that is free to rotate about an axis through its centre and 
is moving in a plane that is perpendicular to its axis of rotation. The 
coordinates of the centre of 
mass of the ball will be denoted by $(x,y)$ and we denote its angular speed as 
$\omega$. It is assumed for simplicity that there is a single point of contact 
between the ball and a compliant 
surface. The contact force acting on the ball can be decomposed 
into the normal and tangential components, denoted by $\lambda_N$ and 
$\lambda_T$ respectively. See Figure \ref{fig:setting_notation} for an  
illustration.

We shall consider different models for the visco-elastic compliance of the 
surface, starting from the simplest models of finite-time impact with friction
\cite{stronge_impacts}, to those that allow for
a finite time interval of contact from impact to lift off, under assumptions 
about the  visco-elastic properties of the surface. To establish some notation, 
we shall choose an origin of co-ordinates so that at the moment of impact, the 
ball is at the origin, so that $(x_0,y_0)=(0,0)$ and is travelling with initial 
horizontal velocity $\dot{x}_0$ and vertical velocity $\dot{y}_0<0$. We shall 
henceforth choose a length scale such that the ball has unit radius $R=1$. 
Similarly, the spin in the plane of motion at impact is denoted by $\omega_0$, 
which can either be positive (counter-clockwise rotation, or 
\textit{backspin}) or negative (clockwise rotation, or \textit{topspin}). 
This quantity will always be expressed in $\mathrm{rad\,s}^{-1}$. In this study 
the ball is modelled as a rigid sphere.

\subsection{Experimental observations}

\begin{figure}[!htb]
	\centering
	\begin{subfigure}{0.45\linewidth}
		\centering
		\includegraphics[width = \linewidth]{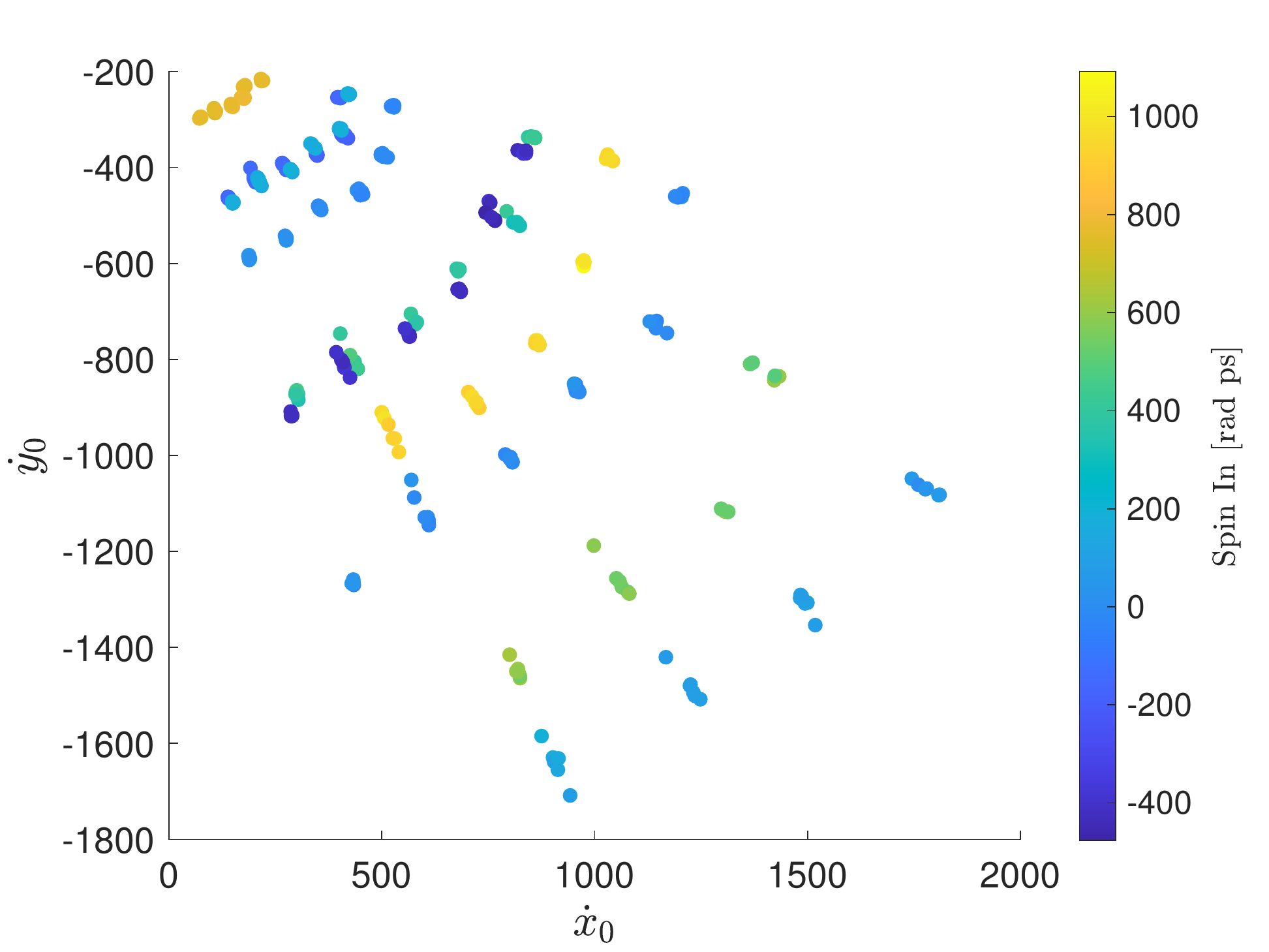}
		\caption{ }
		\label{sfig:campA_ic}
	\end{subfigure}
	\begin{subfigure}{0.45\linewidth}
		\centering
		\includegraphics[width = \linewidth]{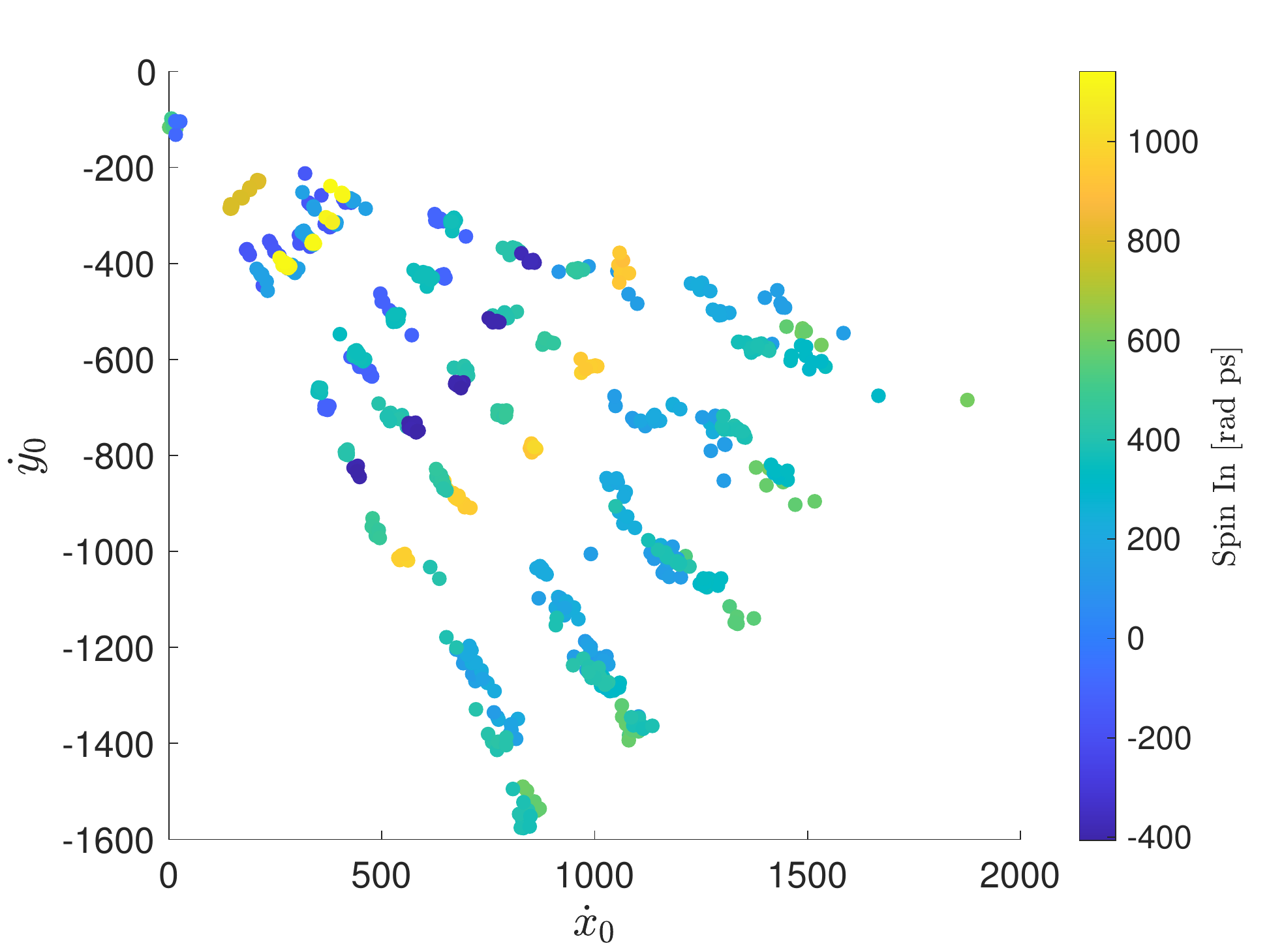}
		\caption{ }
		\label{sfig:campB_ic}
	\end{subfigure}
	\caption{ An overview of the span of the initial condition space. Note that 
		data is dimensionless here ($R=1$). (a) 
		Initial conditions obtained for the testing of artificial turf; (b) 
		Initial conditions obtained for the teeing turf.}
	\label{fig:ics}
\end{figure}

The main motivation for this paper is the recent set of experimental
measurements obtained on golf ball bounce, the details of which are
presented in the companion paper \cite{experiments}. Using a bespoke launcher,
golf balls with a wide variety of incoming velocities and spins were launched 
at two different surfaces: an artificial turf and a well-maintained teeing 
turf. The results were captured using high-speed photography at frame rates of 
up to 10,000 frames per second. Figure~\ref{fig:vid_still} shows an example 
field of view. Automatic image processing enabled accurate measurement of time 
series of velocities and spin within the plane of the ball's motion. In all 
cases it was found that the ball remained in contact with the surface for a 
finite amount of time, during which clear deformation of the ground could be 
seen, by noting that the position of the top of the ball fell and then rose. 
The range of initial conditions covered by the experimental campaign can be 
seen in Table \ref{tab:ics}, and its visual representation is given in 
Figure~\ref{fig:ics}.

\begin{table}[]
	\caption{Range of initial conditions used in the experimental campaigns. 
		For nondimensionalisation of data see the main text. }
	\label{tab:ics}
	\begin{tabular}{|l|ll||ll|}
		\hline
		\multirow{2}{*}{} & 
		\multicolumn{2}{c||}{Astroturf}                      
		&
		\multicolumn{2}{c|}{Real 
			turf}\\ 
		\cline{2-5} 
		& \multicolumn{1}{c|}{Dimensional range}                      & 
		Dimensionless range & \multicolumn{1}{c|}{Dimensional 
			range}                      & Dimensionless range \\ \hline
		$\dot{x}_0$       & \multicolumn{1}{l|}{$(1.53,\,38.6) \, [\mbox{m 
				s}^{-1}]$}   & $(71.1,\,1810)$  & 
				\multicolumn{1}{l|}{$(0.0171,\,36.8) 
			\, [\mbox{m s}^{-1}]$} & $(0.801, \, 1720)$   \\ \hline
		$\dot{y}_0$       & \multicolumn{1}{l|}{$(-36.4,\,-4.61) \, [\mbox{m 
				s}^{-1}]$} & $(-1710,\,-216)$ & 
				\multicolumn{1}{l|}{$(-33.0,\,-2.14) \, 
			[\mbox{m s}^{-1}]$} & $(-1550,\,-101)$     \\ \hline
		$\omega_0$        & \multicolumn{1}{l|}{$(-4550,\, 1040) \, 
			[\mbox{rpm}]$}      & $(-477,\,1090)$  & 
			\multicolumn{1}{l|}{$(-3880,\, 
			1090) \, [\mbox{rpm}]$}      & $(-407,\,1140)$      \\ \hline
	\end{tabular}
\end{table}

\begin{figure}[!htb]
	\centering
	\begin{subfigure}{0.45\linewidth}
		\centering
		\begin{overpic}[width = 0.8\linewidth]{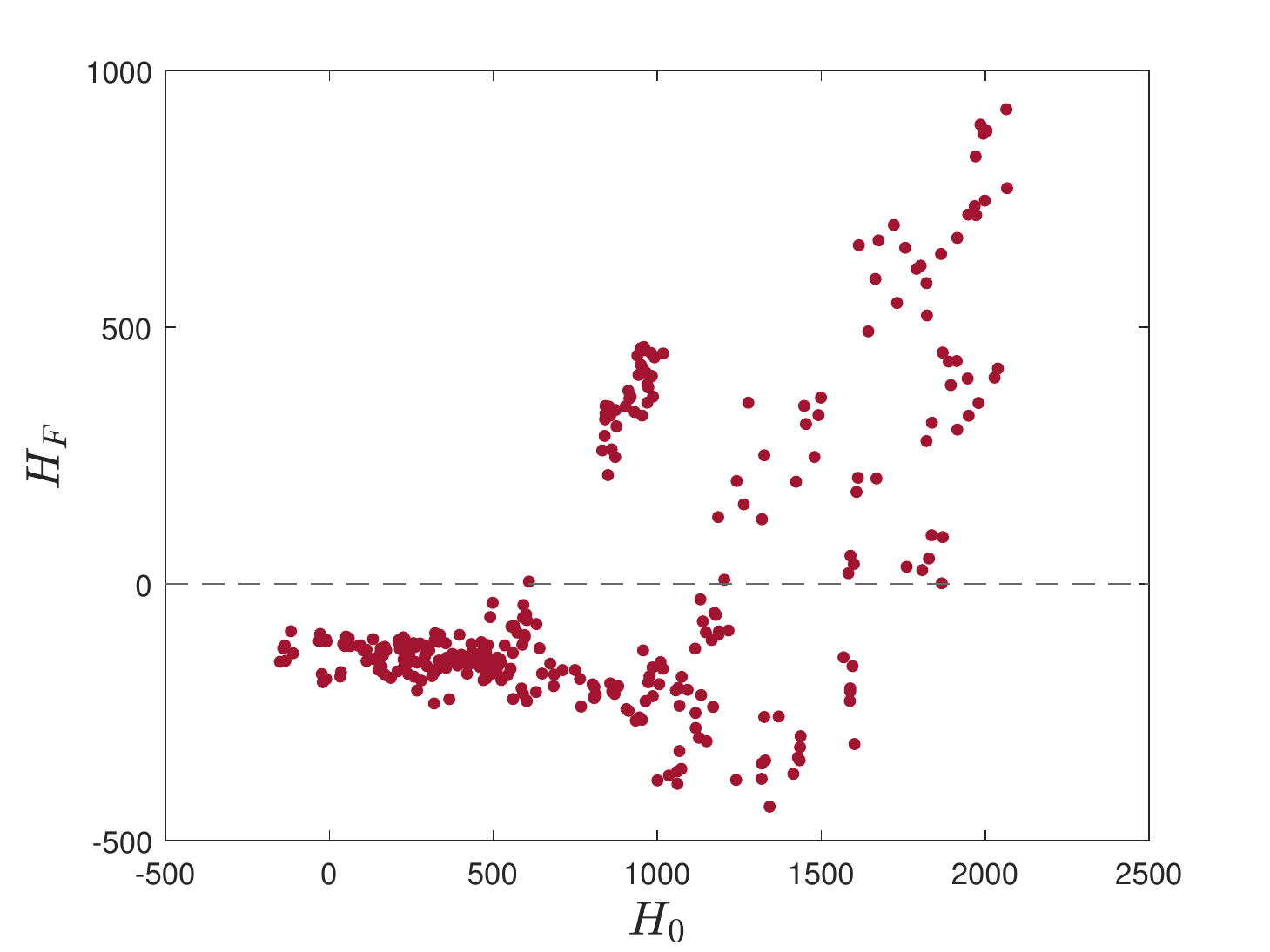}
			
		\end{overpic}
		\caption{}
		\label{sfig:H_data_astro}
	\end{subfigure}
	\begin{subfigure}{0.45\linewidth}
		\centering
		\begin{overpic}[width = 0.8\linewidth]{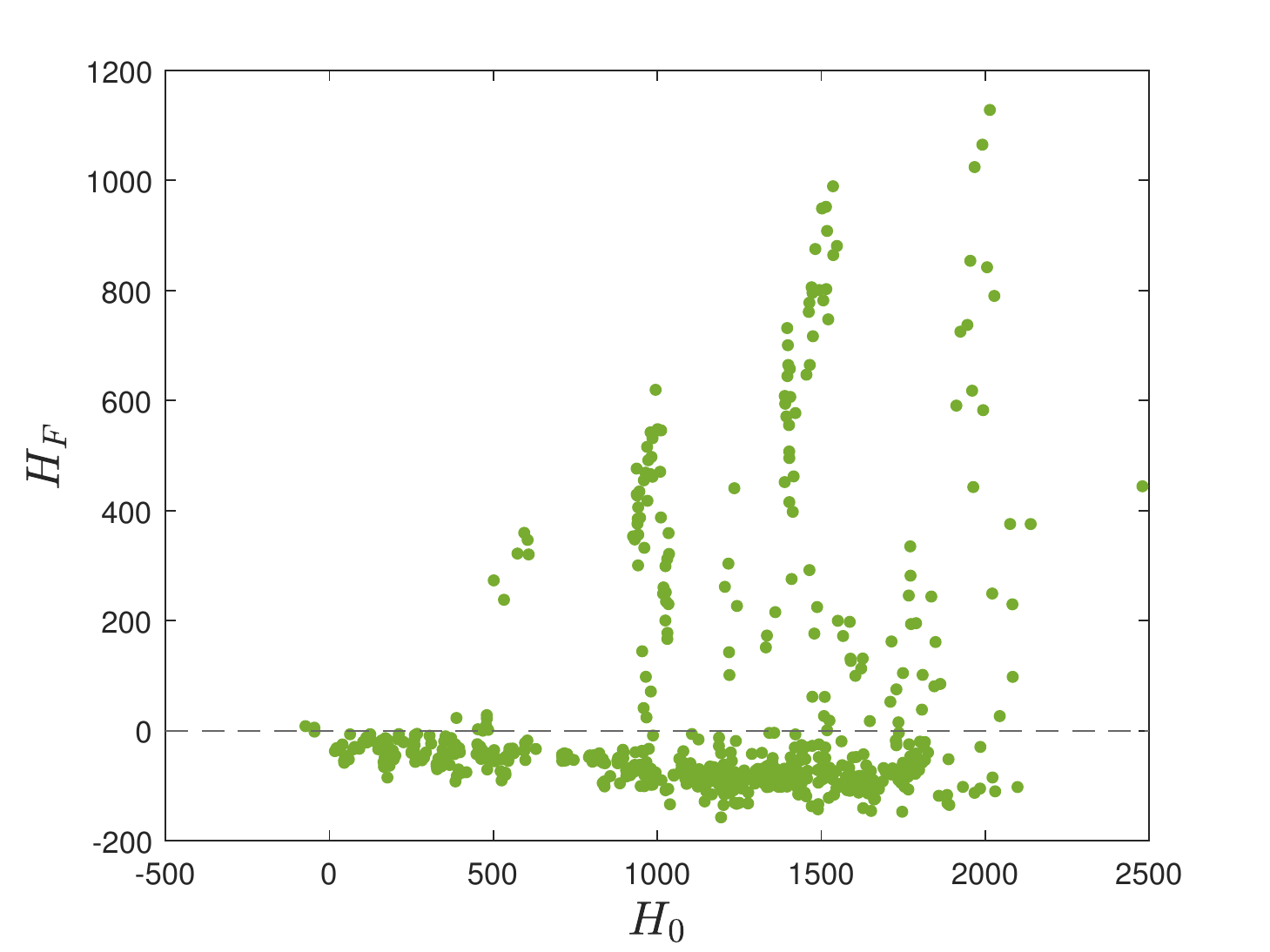}
			
		\end{overpic}
		\caption{}
		\label{sfig:H_data_turf}
	\end{subfigure}
	\caption{Scatter of the experimental results of tangential velocity at 
		lift off versus tangential velocity at first contact; (a) using 
		artificial turf, 
		(b) using fairway turf.}
	\label{fig:tangential_vel_experiments}
\end{figure}

Figure \ref{fig:tangential_vel_experiments} presents just one aspect of the 
data, namely the relation between the incoming and outgoing relative tangential 
velocity $H(t)$ between the ball and the surface at the start and the end of 
the bounce. The results are presented using the notation of 
Sec.~\ref{sec:analysis} below, where (recalling that we rescale length so that 
$R=1$ in our setting) we have that
\begin{equation}
	H(t) = \dot{x}(t) + \omega(t).
\end{equation} 
We also use the notation $H_0 = \dot{x}_0 + \omega_0$ to denote the tangential 
velocity component at the moment of impact and $H_F = \dot{x}_F + \omega_F$ 
at the moment of lift off.

\begin{figure}[!htb]
	\centering
	\begin{subfigure}{0.45\linewidth}
		\centering
		\begin{overpic}[width = \linewidth]{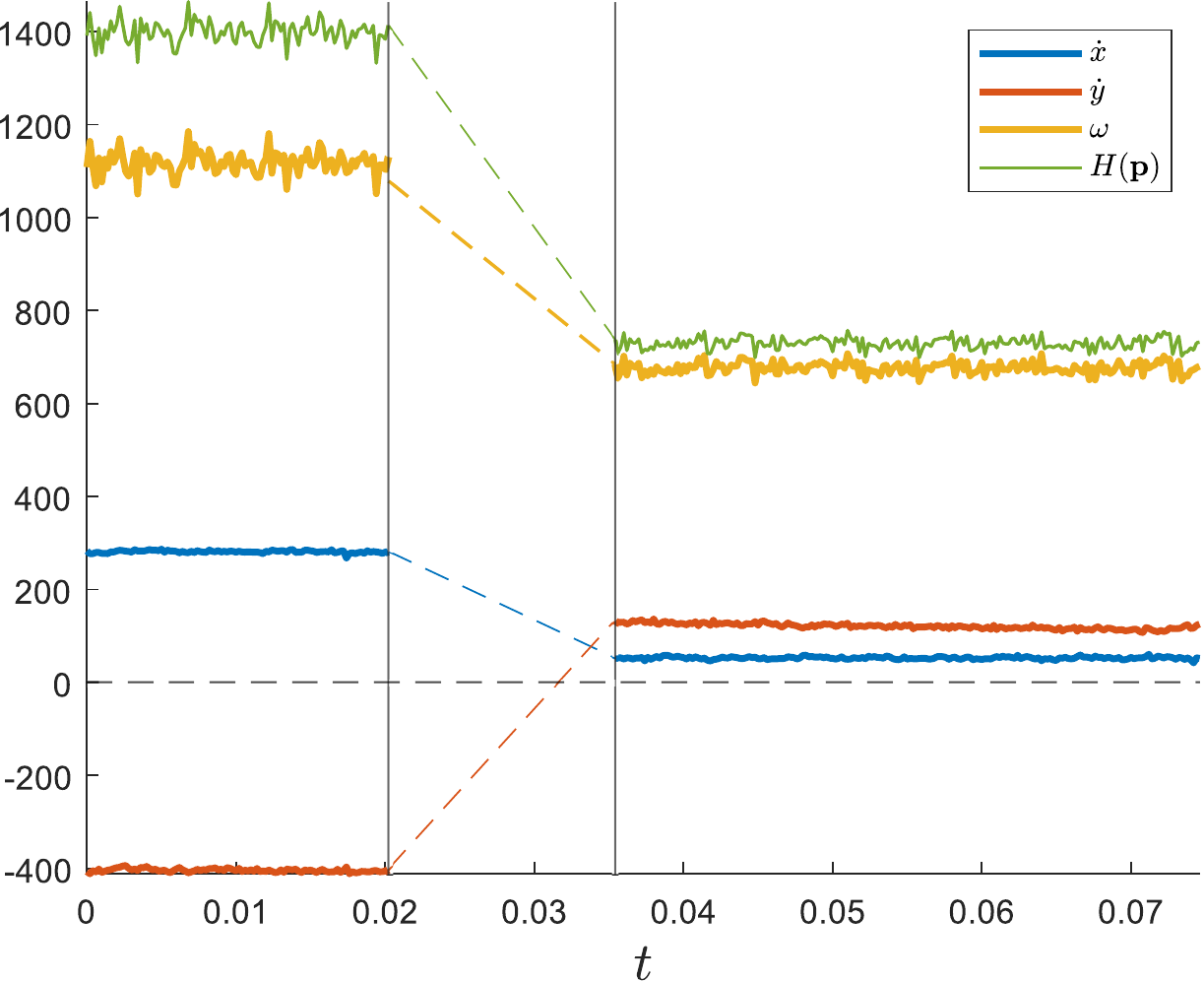}
			\put(33.5,40){\tiny Bounce}
		\end{overpic}
		\caption{}
		\label{sfig:data_spin_high}
	\end{subfigure}
	\begin{subfigure}{0.45\linewidth}
		\centering
		\begin{overpic}[width = \linewidth]{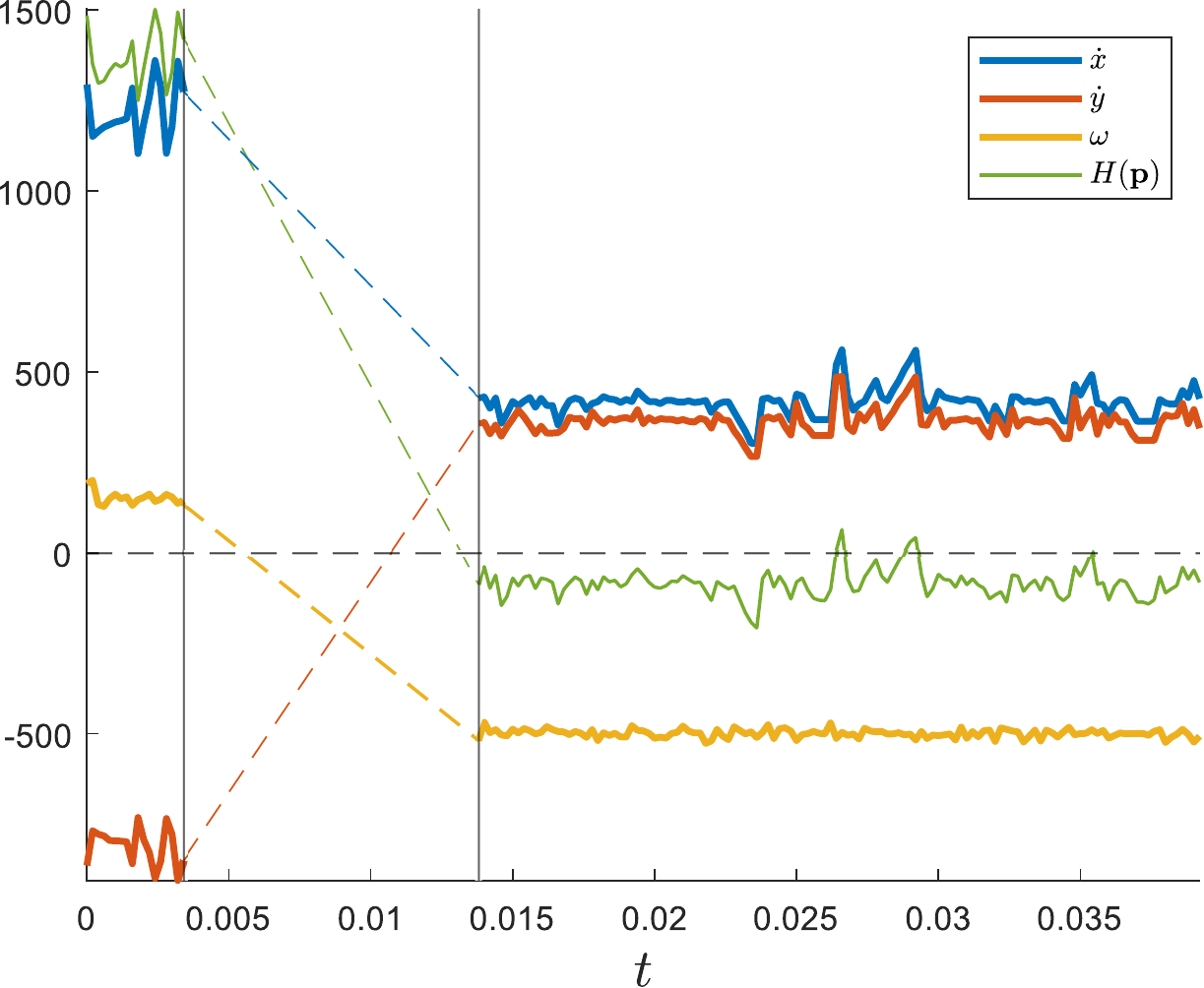}
			\put(20,40){\tiny Bounce}
		\end{overpic}
		\caption{}
		\label{sfig:data_spin_low}
	\end{subfigure}
	\caption{Changes in the (non-dimensional) velocities before and after the 
		bounce. 
		(a) High tangential velocity at	lift off. (b) Negative 
		tangential velocity at lift off. All data gathered from the 
		experimental session using teeing turf. Colours blue, red, orange and 
		green represent	vertical, horizontal, angular and tangential 
		velocities respectively.}
	\label{fig:data_vel_change}
\end{figure}

We observe that for both the synthetic and 
actual turf, the data falls roughly into two classes. For one class, the 
tangential velocity at lift off is large and positive. Further examination of 
the data reveal that these cases correspond to trajectories that land with a 
large amount of backspin (that is, a large positive angular velocity), slip 
throughout the bounce impact, and lift off with a reduced
amount of backspin. Another somewhat larger class of samples, appear to lift 
off with small negative relative velocity. Analysis of the trajectories show 
that these trajectories typically enter with a small amount of backspin, no 
spin or topspin but lift off with a small amount of topspin. In particular, 
there is little evidence of balls lifting off rolling, that is with the 
relative tangential velocity being zero. Examples of the raw time-series data 
for each of the two cases are illustrated in  Figure \ref{fig:data_vel_change}. 
Note how the changes during the bounce are not recorded, as the ball is not 
fully visible in that phase due its immersion within the surface.

\subsection{Rigid impact with friction} 
\label{sec:rigid_bounce}

Let us first look at the case of a purely rigid bounce, that is where the rigid 
ball bounces off a rigid surface subject to Coulomb friction. During the bounce 
normal force due to the contact, gravity and tangential force due to Coulomb 
friction act on the body, yielding the equations of motion
\begin{equation}
	\ddot{x} = \lambda_T, \qquad \ddot{y} = \lambda_N - g, \qquad \dot{\omega} 
	= \frac{5}{2} \lambda_T.
\end{equation}
Many authors have considered impact problems of this type and results typically
only differ by the assumptions made about how to capture energy loss through
a coefficient of friction; see e.g.~\cite{Brogliato} and references therein. 

Here we follow the approach used in \cite{stronge_impacts} where the impact is 
assumed to last an infinitesimal amount 
of time and the contact forces are assumed to be $\delta$-function-like 
distributions. An important part of the process is to model transitions between 
roll and slip during the infinitesimal time of impact. This can
lead to subtle effects if the impacting body is slender \cite{danko1}.  In the 
case of sphere, however, the implementation of an analogue to the single 
degree-of-freedom coefficient of restitution is straightforward and
most methods lead to the same answer \cite{Battle}. In particular here  we 
assume, given the landing normal velocity $\dot{y}_0$, that lift-off is deemed 
to occur when the ball reaches a  normal velocity 
$\dot{y}_F = -r \, \dot{y}_0$, where $r$ is called a kinematic coefficient of 
restitution. For simplicity, we assume $r$ is constant (although note
the empirical studies \cite{penner2002,Roh} which suggests golf turf may be 
modelled by a velocity-dependent coefficient of restitution). We are concerned 
about the dynamics of the point of contact, at which we can write the 
tangential and normal velocities as 
\begin{equation}
	\vec v = [v_Y, v_N]^{\intercal} = [\dot{x} + \omega, \dot{y}]^{\intercal},
\end{equation}
so that 
\begin{equation}
	\frac{\dd \vec v}{\dd t} =  \begin{bmatrix} \frac{7}{2} & 0 \\ 0  & 1 
	\end{bmatrix} \begin{bmatrix} \lambda_T \\ \lambda_N \end{bmatrix} = m^{-1} 
	\vec \lambda,
\end{equation}
where $\vec \lambda$ is the vector of contact forces and $m^{-1}$ can be 
thought of as the inverse of the mass matrix in the ``coordinate system'' of 
$(v_T, v_N)$.

Let us denote by $t=t_0$ the time at which the impact is initiated. The 
subscript $0$ will denote evaluation of other quantities at $t_0$. We consider 
small time interval during which the ball is in contact with the rigid surface 
such that for all $t$ in that interval $t-t_0= \ve \tilde{t}$, where $\ve \ll 
1$. During that time we assume $\ve \vec \lambda, v_T, v_N = \Ord{1}$, and 
denote $\tilde{ \vec\lambda} = \ve \vec \lambda$, so that 
\begin{equation}
	\frac{\dd \vec v}{\dd t} = m^{-1} \tilde{\vec \lambda} + \Ord{\ve}.
\end{equation}
We now follow the idea introduced by Stronge in \cite{stronge_impacts}, where 
we replace our independent time variable by the new one, namely \textit{the 
	normal impulse}
\begin{equation}
	n = \int \lambda_N \, \dd \tilde{t},
\end{equation}
which is a strictly increasing function of $\tilde{t}$. Dropping the tilde (for 
simplicity of notation) and ignoring terms $\Ord{\ve}$ we now have, with 
respect to the new independent variable,
\begin{equation}\label{eq:v_divs_gen} 
	\frac{\dd \vec v}{\dd n}= \frac{1}{\lambda_N} m^{-1} \vec \lambda = 
	\begin{bmatrix} \frac{7}{2} \frac{\lambda_T}{\lambda_N}\\
		1 \end{bmatrix}.
\end{equation}

During the slipping phase of the motion $\lambda_T  = - \sign{ v_T }  \mu
\lambda_N $, where $\mu$ is the coefficient of friction. 
Substituting this into Equation \eqref{eq:v_divs_gen} and rearranging the 
derivatives to eliminate $n$ we have that during the slipping phase of the 
rigid bounce
\begin{equation}\label{eq:v_div_slip}
	\frac{\dd v_T}{\dd v_N} = \pm \frac{7}{2}\mu.
\end{equation}

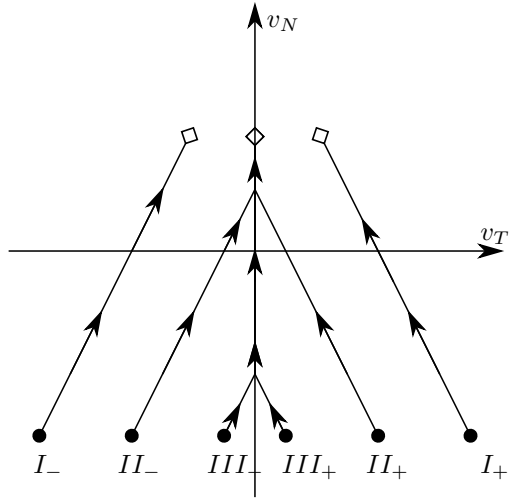
\begin{figure}[!htb]
	\centering
		\begin{overpic}[width=0.4\textwidth]{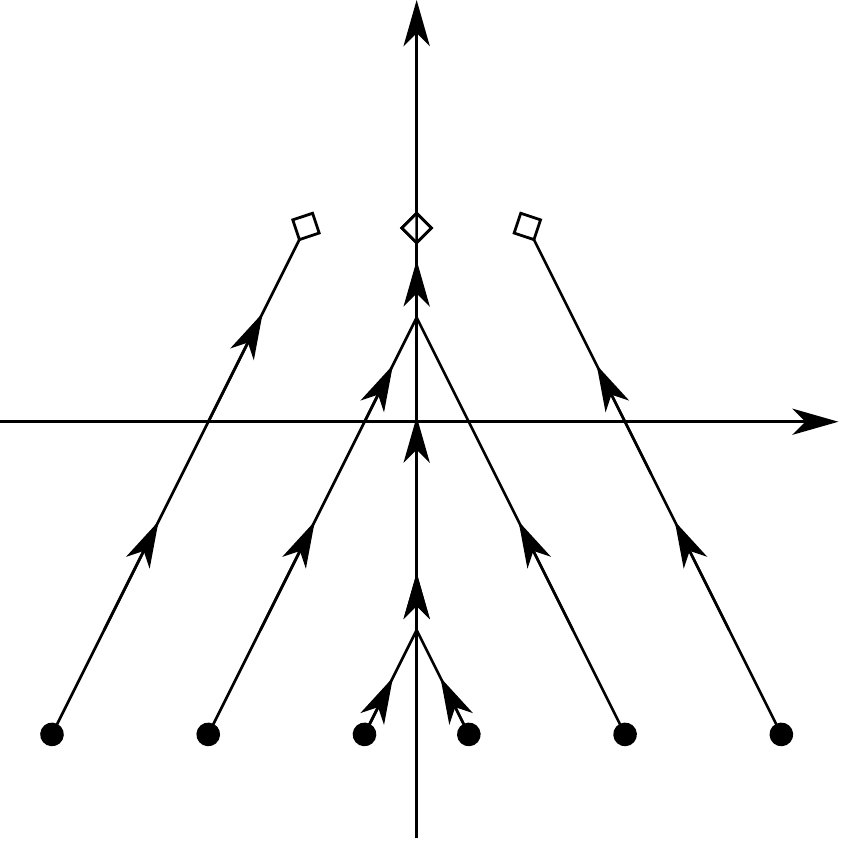}
			\put(95,52){$v_T$}
			\put(52,95){$v_N$}
			\put(5,5){$I_{-}$}
			\put(22,5){$II_{-}$}
			\put(40,5){$III_{-}$}
			\put(95,5){$I_{+}$}
			\put(72,5){$II_{+}$}
			\put(55,5){$III_{+}$}
		\end{overpic}

	\caption{A graphic interpretation of changes in tangential and normal 
		velocities during impact (Equations \eqref{eq:v_div_slip} and 
		\eqref{eq:v_div_roll}). The start of the impact is denoted with solid 
		black 
		dot, and the end is denoted with a white rhombus. We differentiate 
		between 
		cases of slip throughout the impact ($I_{\pm}$) and slip followed by 
		stick 
		in the restitution ($II_{\pm}$) or compression ($III_{\pm}$) phase.}
	\label{fig:trajectories}
\end{figure}

However, defining the tangential force during the roll needs a bit more care. 
We note that during the roll $v_T = \dot{x} + \omega =0$ and thus it must also 
be true that $\ddot{x} + \dot{\omega} = \lambda_T+\frac{5}{2} \lambda_T = 
\frac{7}{2} \lambda_T=0.$ This in turns requires $\lambda_T=0$ during the roll, 
which implies that once the ball enters rolling, it may never leave it, and 
during the roll. Combined with Equation \eqref{eq:v_divs_gen} we have the for 
the rolling stage of the bounce
\begin{equation}\label{eq:v_div_roll}
	\frac{\dd v_T}{\dd v_N} = 0.
\end{equation}
With these equations it is straight forward to work out conditions under which 
$v_N$ and $v_T$ vary, which can be expressed graphically using the method 
presented in
\cite{danko1}, see in Figure 
\ref{fig:trajectories}. In particular, we distinguish between three cases
in the space of initial conditions 
$(\dot{x}_0,\dot{y}_0,{\omega}_0)$:

\noindent \textbf{Case I:} the ball will slip through the impact (case 
$I_{\pm}$) if 
\begin{equation*}
	\left| \frac{\dot{x}_0 + \omega_0}{\dot{y}_0}\right| > \frac{7}{2}\mu 
	(1+r);
\end{equation*}

\noindent \textbf{Case II:} the ball will enter the roll in the restitution 
phase (case 
$II_{\pm}$) if 
\begin{equation*}
	\frac{7}{2}\mu < \left| \frac{\dot{x}_0 + \omega_0}{\dot{y}_0}\right| < 
	\frac{7}{2}\mu (1+r);
\end{equation*}

\noindent \textbf{Case III:} the ball will enter the roll in the compression 
phase (case 
$III_{\pm}$) if 
\begin{equation*}
	\left| \frac{\dot{x}_0 + \omega_0}{\dot{y}_0}\right| < \frac{7}{2}\mu.
\end{equation*}

Note that none of the cases allows the ball to enter slipping once it has begun 
to roll during the bounce.  This contravenes what we saw in the experimental 
data. Specifically, if we consider balls entering bounce with $v_T>0$ then, 
this theory allows lift-off from bounce to occur either in forward slip 
($v_T>0$) or in roll $(v_T=0)$; there is no initial condition that leads to 
`spin reversal' during bounce. (Note though that this observation that this 
does not actually preclude a `backwards bounce' because with high backspin a 
ball that has $v_T=0$  would lift off with $\dot{x}<0$).

Note that the absence of spin reversal also occurs if we introduce normal 
compliance into the theory, so that the bounce occupies a short 
$\Ord{\ve}$ period of time, while keeping the Coulomb friction 
assumption in the tangential direction. For example, a straightforward 
calculation (see \cite{thesis_prep} for details) can be performed 
assuming a Kelvin-Voigt-like model
of a spring and damper connected in parallel.  Assuming the stiffness of the 
spring to be $1/\ve$ and 
the damping parameter of the dashpot to be $d/\ve$, where $d$ is 
a damping ratio, then one obtains the same 
conclusion about the absence of spin reversal. In particular, in the limit 
$\varepsilon \to 0$, with $d< \frac{1}{2}$ (to assure the underdamped 
solution), one recovers the rigid bounce theory, with the 
coefficient of restitution 
\begin{equation}
	r \sim \exp \left[ -\frac{\pi}{2} d + \frac{1}{2} d^2 + \Ord{d^3} \right] + 
	\Ord{\ve}.
\end{equation}
More generally, it is clear that in order to observe spin reversal in a simple 
point-contact model, a more complex model of the mechanics in the tangential 
degree of freedom is required.

\section{Linear elasto-plastic model}
\label{sec:two_cases}

\begin{figure}[!htb]
	\centering
	\begin{overpic}[width = 0.4\linewidth]{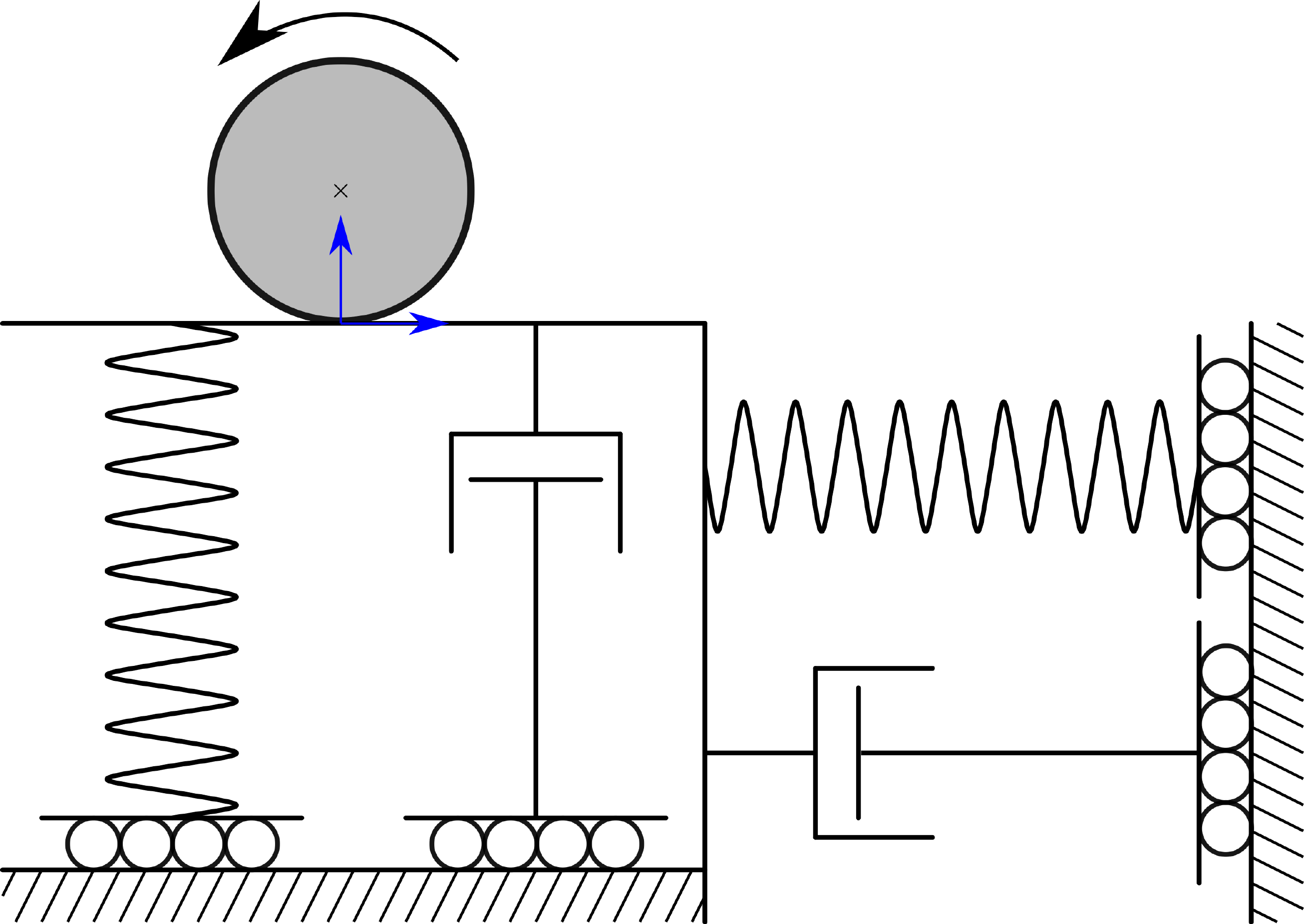}
		\put(20,58){\tiny $(x, y)$}
		\put(30,71){\tiny $\omega$}
		\put(42,20){\tiny $\frac{2d_2}{\ve_2}$}
		\put(2,20){\tiny $\frac{1}{\ve_2^2}$}
		\put(73,6){\tiny $\frac{2d_1}{\ve_1}$}
		\put(73,46){\tiny $\frac{1}{\ve_1^2}$}
		\put(19,51){\tiny $\lambda_N$}
		\put(34,48){\tiny $\lambda_T$}
	\end{overpic}
	\caption{A linear elasto-plastic bounce mechanism based on the Kelvin-Voigt 
		model}
	\label{fig:settings_kv}
\end{figure}
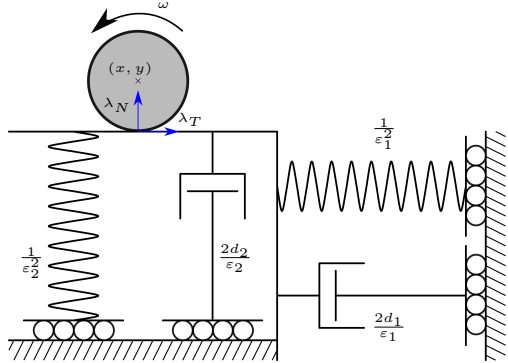

Let us now consider a model where the ball remains rigid as before, but 
compliance in both horizontal and vertical directions can be observed. 
As depicted in Figure \ref{fig:settings_kv}, 
the elasto-plastic behaviour is assumed to follow one of a spring and 
damper connected in parallel, known in the literature as the Kelvin-Voigt 
model, see e.g.~\cite{findley_creep}. The dynamics between the normal 
and vertical directions is decoupled, and it is assumed that
the ball is in contact with the 
surface at a single point only, where the normal and tangential forces again 
follow the Coulomb friction law. The equations of motion for such a system can 
be written in dimensionless form 
\begin{equation}\label{eq:kv}
	\begin{aligned}
		\ddot{x} + \frac{2d_1}{\ve_1} \dot{x} + \frac{1}{\ve_1^2} x & = 
		\lambda_T,\\
		\ddot{y} + \frac{2d_2}{\ve_2} \dot{y} + \frac{1}{\ve_2^2} y & = 
		-g,\\
		\dot{\omega} & = \frac{5}{2}\lambda_T.
	\end{aligned}
\end{equation}
Here $g$ is the acceleration due to gravity, $\ve_{1,2}$ are the stiffness 
ratios in horizontal and vertical directions and $d_{1,2}$ are the damping 
ratios in the respective directions. Note the scaling that $\ve_{1,2} \ll 1$ is 
consistent with a time scale in which the bounce take places over a rapid time 
scale, and passing to the 
limit $\ve_{1,2} \to 0$ will result in a rigid model with a horizontal and 
vertical coefficients of restitution (determined by the damping ratios 
$d_{1,2}$, provided that $d_{1,2}<1$ so that the system is underdamped. In 
particular, we will assume $d_2<1$, which 
will insure  that a ball coming into contact with the ground with 
$\dot{y}_0<0$ will lift-off at some later time with $\dot{y}_F>0$. 

Specifically, we can readily compute that the ball 
will remain in contact with the compliant surface for as long as 
\begin{equation}
	\lambda_N = -\frac{2d_2}{\ve_2} \dot{y} - \frac{1}{\ve_2^2} y -g>0
\end{equation}
and will lift off when 
$\lambda_N=0$ 

Returning to the physical motivation of golf-ball turf interaction, in addition 
to the scaling that $\ve_{1,2} \ll 1$, we shall additionally assume that $\ve_1 
\ll \ve_2$. Such a scaling implies that there is far greater compliance in the 
normal direction than the tangential. That  this is a natural assumption can be 
understood by briefly considering instead the motion that would occur if 
$\ve_1 = \ve_2$, so that the stiffness experience is isotropic. In that case, 
in the absence of spin or damping, and assuming $\dot{x}$ a ball with inbound 
angle $\phi = \arctan (\dot{x},\dot{y})$ would move along a straight line in 
the $(x,y)$-plane and would lift-off 
with angle $\phi+\pi$.  This is not how turf behaves. Instead, observational 
data suggest that such a ball would lift with outbound angle close to $-\phi$. 

With such scaling, it is useful to introduce a new time rescaling 
\begin{equation}
	\tau = \ve_2^{-1} t = \Ord{1}
\end{equation}
and a dimensionless parameter 
\begin{equation}
	\eta = \ve_2 / \ve_1 \ll 1.
\end{equation}

The equations of motion can now be written in the dimensionless form 
\begin{equation}\label{eq:kv_rescaled}
	\begin{aligned}
		X'' + 2d_1 \eta X' + \eta^2 X & = \Lambda_T,\\
		Y'' + 2d_2 Y' + Y & = - \ve_2^2 g,\\
		\Omega' & = \frac{5}{2} \Lambda_T,
	\end{aligned}
\end{equation}
where $X(\tau) = x(t)$, $Y(\tau) = y(t)$, $\Omega(\tau) = \omega(\tau)$, 
$\Lambda_T=\eps_2^{-2} \lambda_T$, and a prime represents differentiation with 
respect to $\tau$. We also define 
$\Lambda_N = \ve_2^{-2}\lambda_N,$ that is
\begin{equation}
	\Lambda_N = - 2d_2 Y' - Y -\ve_2^2 g>0.
	\label{eq:nforce}
\end{equation}

The Coulomb friction law (in the new variables) dictates that for
$\lvert  \Lambda_T \rvert = \mu \Lambda_N$ the ball slips on the surface. 
However,  
greater care is needed for the rolling motion, that is when $\lvert \Lambda_T 
\rvert = \mu \Lambda_N$. 

\subsection{Defining rolling; a Filippov formulation} \label{sec:filippov}

To consider what happens during rolling, 
it is useful to model the system as a Filippov System 
\cite{piecewisebook,filippov,jeffrey_hidden}. That is, with $\Lambda_T$ well 
defined for the cases of slip, we note that the equations 
\eqref{eq:kv_rescaled} can be written in the form
\begin{equation}
	\p ' = 
	\begin{cases}
		F_1(\p) & \mathrm{if } \quad H(\p) > 0,\\
		F_2(\p) & \mathrm{if } \quad H(\p) < 0,
	\end{cases}
\end{equation}
where $\p$ is the vector of state variables. Specifically, we
have 
\begin{equation} \label{eq:vec_p}
	\p = [X, X', Y, Y', \Omega]^{\intercal}
\end{equation}
and the smooth function
\begin{equation}
	H(\p) = X' + \Omega.
\end{equation}
In 1988 Filippov \cite{filippov} proposed a consistent method 
of defining dynamics along the surface of discontinuity 
$\{H(\p)=0\}$ in which a new vector field is introduced
\begin{equation}\label{eq:sliding}
	F_s(\p) = (1-\ga) F_1(\p) + \ga F_2 (\p) \quad \mathrm{ when } \quad H(\p) 
	= 0,
\end{equation}
where $\ga \in [0,1]$. That is, the vector field along the discontinuity  
surface 
is the unique convex linear combination of the two  neighbouring vector fields, 
such that the resulting vector lies along the  discontinuity. The vector field 
$F_s$ in our system would define rolling, but confusingly, in deference to 
control applications, in Filippov systems $F_s$ is termed the {\em sliding} 
vector field.

A `sliding' trajectory is allowed to leave the surface of  discontinuity when 
either $\ga = 0$ (therefore leaving into the region $H(\p)>0$ where vector 
field $F_1$ applies) or when $\ga =1$ (where $F_2$ applies). At the same time, 
$\ga$  can be explicitly calculated during the sliding motion by making that 
the assumption that the vector field must be orthogonal to $\nabla H$. Thus, we 
have
\begin{equation}\label{eq:alpha}
	\ga(\p) = \frac{F_1(\p) \cdot \nabla H}{(F_1(\p) - F_2(\p) )
		\cdot \nabla H}.
\end{equation}

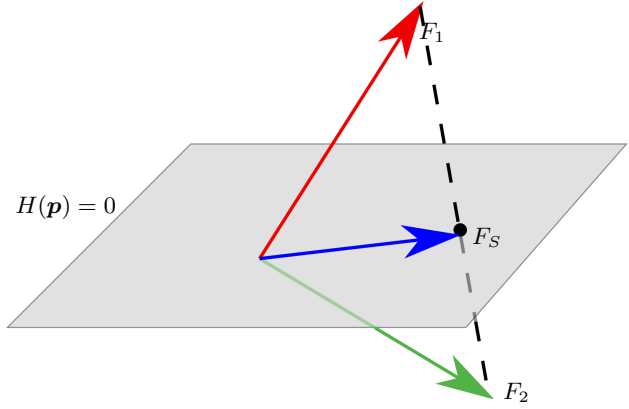
\begin{figure}
	\centering
	\begin{overpic}[width=.5\textwidth]{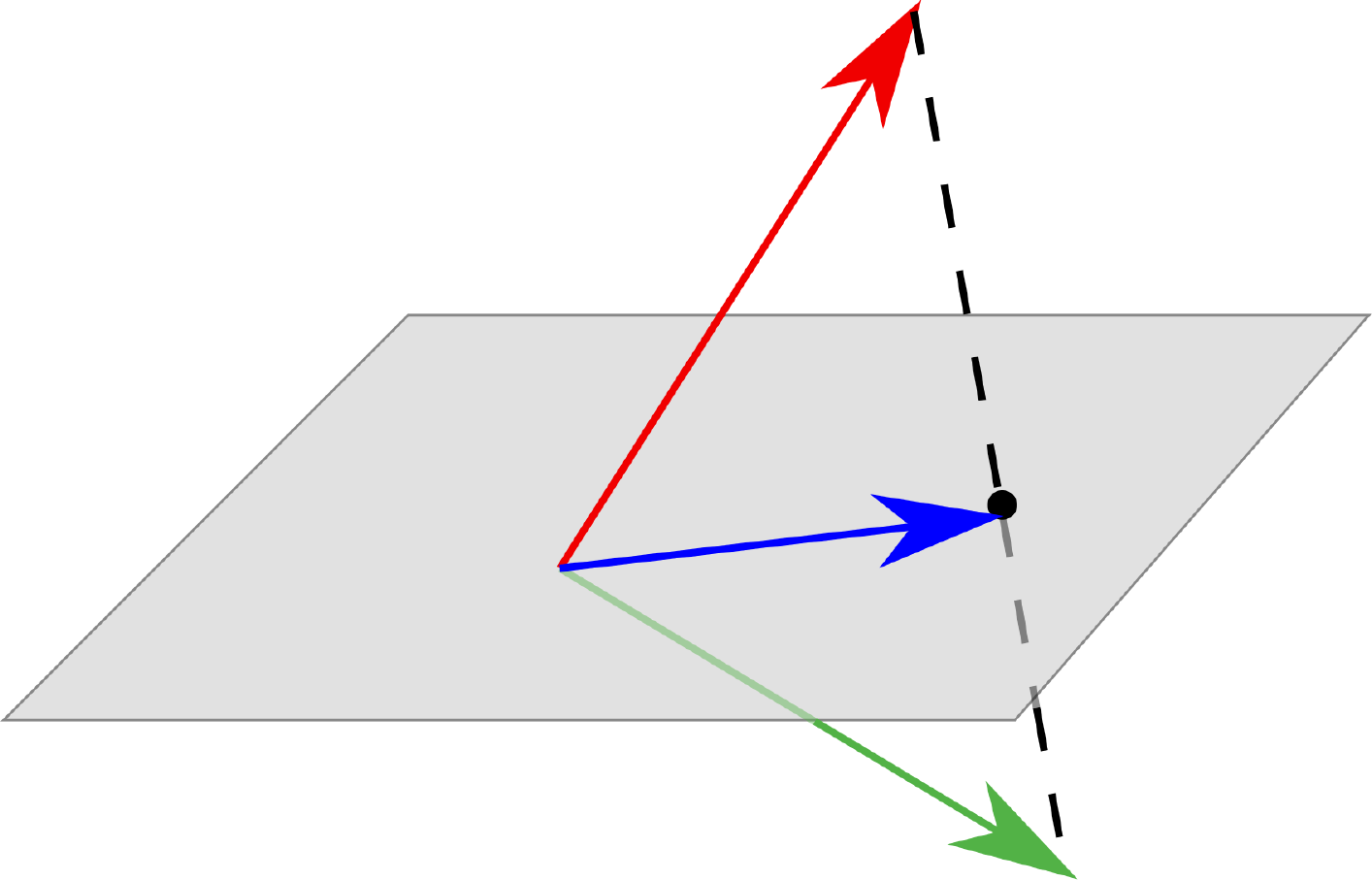}
		\put(80,0){\small $F_2$}
		\put(65,58){ \small $F_1$}
		\put(0,30){ \small $H(\p) =0 $}
		\put(75,25){\small $F_S$}
	\end{overpic}
	\caption{Definition of the sliding dynamics --- the  discontinuity surface 
		is denoted here with a plane ($H(\p) = 0$), the  dynamics according to 
		the 
		two vector fields $F_1$ and $F_2$ are denoted  with red and green 
		arrows 
		(above and below the discontinuity surface).  Their convex linear 
		combination lies along the dashed line, and we chose the unique 
		solution 
		that also aligns with the sliding surface.}
	\label{fig:filippov_convexity}
\end{figure}

\subsection{Analysis of roll, slip and lift-off transitions}

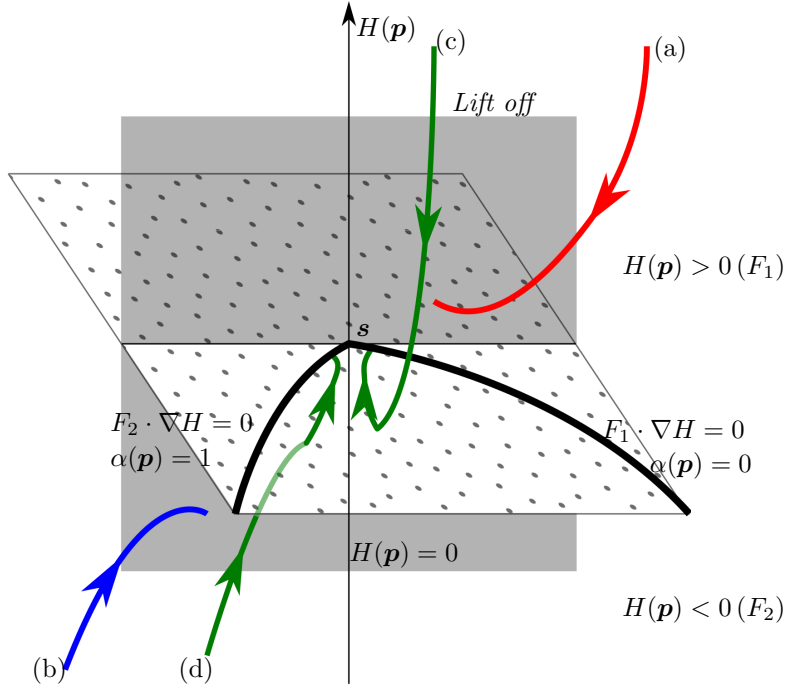
\begin{figure}
	\centering
	\begin{overpic}[width = 0.55\linewidth]{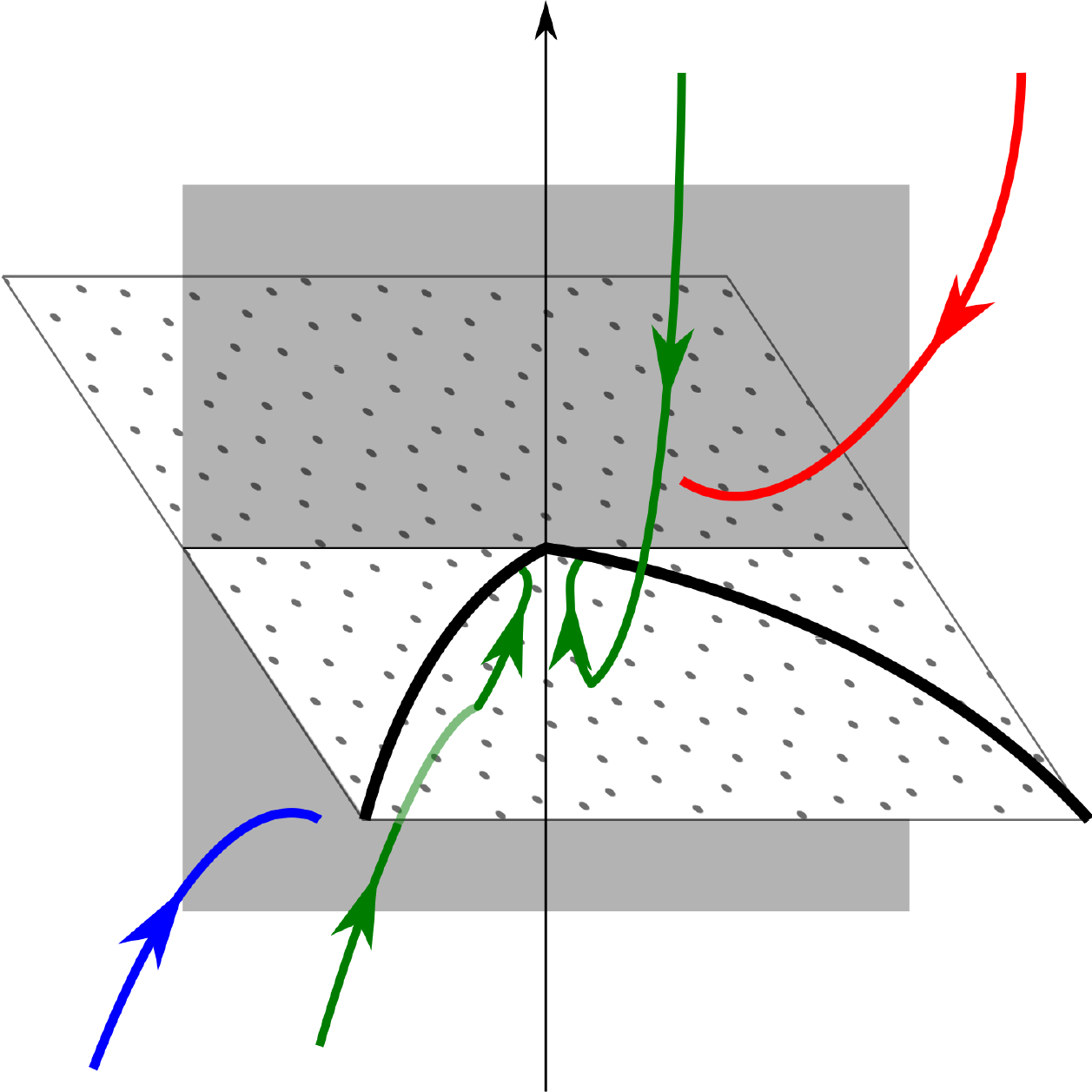}
		\put(51,95){$H(\p)$}
		\put(50,18){$H(\p) = 0 $}
		\put(65,84){{\it Lift off}}
		\put(90,60){$H(\p) > 0 \, (F_1) $}
		\put(90,10){$H(\p) < 0 \, (F_2) $}
		\put(15,37){$F_2 \cdot \nabla H = 0 $}
		\put(15,32){$\ga (\vec p) = 1 $}
		\put(87,36){$F_1 \cdot \nabla H = 0 $}
		\put(94,31){$\ga (\vec p) = 0 $}
		\put(94.5,92){(a)}
		\put(3.5,1){(b)}
		\put(62.5,93){(c)}
		\put(25,1){(d)}
		\put(51,51){$\s$}
	\end{overpic}
	\caption{Possible trajectories during the generic elasto-plastic bounce (in 
		reduced geometry). The grey plane represents the lift off condition, 
		the 
		dotted plane is the codimension-1 discontinuity, the solid black curves 
		represent surfaces where $F_1 \cdot \nabla = 0$ and $F_2\cdot \nabla H 
		= 
		0$, which corresponds to the surfaces where a trajectory leaves the 
		surface 
		of discontinuity and enters either $F_1$ or $F_2$ (respectively). Point 
		$\s$ is an intersection of the two switching surfaces (known as a 
		two-fold bifurcation point \cite{jeffrey_update}) and lies on the lift 
		off 
		plane as well.\\
		Red trajectory (a): The ball goes through the bounce with positive slip 
		only ($H(\p)>0$). The bounce is terminated by the intersection with the 
		grey plane, which represents the lift off condition.\\
		Blue trajectory (b): The ball goes through the bounce with negative 
		slip 
		only ($H(\p)<0$). The bounce is terminated by the intersection with the 
		grey plane, which represents the lift off condition.\\
		Green trajectories (c-d): The ball enters the bounce with positive or 
		negative slip (respectively) and eventually reaches the system 
		discontinuity by entering roll ($H(\p) =0$). The bounce dynamics 
		continue 
		along the surface of discontinuity and are to be determined in the next 
		sections. }
	\label{fig:geometries}
\end{figure}

Having now defined the dynamics fully to be 
\begin{equation}
	\p ' = 
	\begin{cases}
		F_1(\p) & \mathrm{if } \quad H(\p) > 0,\\
		F_s(\p) & \mathrm{if } \quad H(\p) = 0,\\
		F_2(\p) & \mathrm{if } \quad H(\p) < 0,
	\end{cases}
\end{equation} 
with $F_S$ specified in Equation \eqref{eq:sliding}, we can 
consider different cases of dynamics during the bounce. We represent possible 
geometries in Figure \ref{fig:geometries}. In it we note that a ball can either 
slip throughout the bounce (either {\em forward slip} 
with $H(\p)>0$ or {\em negative slip} with $H(\p)<0$) 
or can at some time $t$ begin to roll (that is, the dynamics evolves along the 
discontinuity set where $H(\p) =0$).

Let us now direct our attention to the latter of the presented cases, that is 
when the ball is rolling and in particular, we shall study how exiting such 
state is possible. As discussed before, exiting will be dependent on the 
value of the function $\alpha(\p)$, which for our system is explicitly 
calculated to be
\begin{equation}
	\alpha(\p) = \frac{7\mu\,(2 d_2 Y' +Y) - 4d_1 \eta X' - \eta^2 X}{7\mu\,(2 
		d_2 Y'+Y)} + \Ord{\ve_2^2}.
\end{equation}
We now analyse the transitions between rolling and slipping. 

Suppose, for definiteness, that a trajectory along the discontinuity surface
$H(\p) =0$ reaches a point $\p_1$ where $F_1(\p_1) \cdot \nabla H = 
0$. That means that 
the sliding field definition is equal to the vector field $F_1$ at
this point. Moreover, 
as dictated by \eqref{eq:alpha}, the trajectory should leave the surface into 
the half-space
$H(\p_1)>0$. Such a trajectory leaving the discontinuity surface must be 
tangent to 
$H(\p_1) = 0$, therefore we have 
\begin{equation}\label{eq:switching_con}
	F_1(\p_1) \cdot \nabla H = -2 d_1 \eta X' - \eta^2 \, X + \frac{7}{2}\mu 
	\left[2d_2 Y' + Y + \ve_2^2 g\right] =0.
\end{equation}
Notice that the term $2d_2 Y' + Y + \ve_2^2 g $ is in fact the 
$-\Lambda_N$, given by \eqref{eq:nforce} and thus remain negative 
throughout the bounce. On the other hand, in order to balance the first two 
terms on the right-hand side of \eqref{eq:switching_con}, at $\p_1$ we must 
have that  $\Lambda_N = \Ord{\eta} \ll 1$. 

Given this information, let us now consider the possible 
signs of the state variables at $\p_1$. The bounce is initiated at $Y=0$ and 
$Y\leq0$ throughout the bounce. Suppose that at $Y'<0$ at $\p_1$, so that we 
are in the downwards (compression) phase of the bounce.  Then since $\Lambda_N 
= \Ord{\eta}$ and $Y<0$ it must follow 
that both $Y, Y' = \Ord{\eta}$ at $\p_1$. $Y$ being small means we must be 
close to the time of initiation of bounce. However,  initially $Y'(0) = 
\Ord{1}$ which gives a 
contradiction. Thus, we must have that $Y'>0$ at $\p_1$. That is, the 
transition from roll into slip can only happen during the upward (restitution) 
phase of the bounce.

Next, let is compute the second directional derivative at such a point $\p_1$. 
We obtain  
\begin{equation}\label{eq:lie2_kv}
	\begin{aligned}
		\left( F_1 (\p_1)\cdot \nabla \right)^2 H  &= \frac{7}{2}\mu Y' + (7d_2 
		+ 2d_1 \eta) \left[- 2d_2 Y' - Y -\ve_2^2 g\right] + \Ord{\eta^2} \\
		& = \frac{7}{2}\mu 	Y' + (7d_2+2d_1\eta) \Lambda_N + \Ord{\eta^2}.
	\end{aligned}
\end{equation}
Now, the above argument shows that both of the terms on the right-hand side of 
\eqref{eq:lie2_kv} are positive. 
Thus $\left( F_1 (\p_1) \cdot \nabla 
\right)^2 H  >0$, which confirms that the ball will leave the rolling region 
consistently into the $F_1$ region. 

We can similarly analyse we the case where the ball transitions from rolling to 
slipping via the sliding flow becoming tangent to the vector field $F_2$ at 
some point $\p=\p_2$. We now have 
\begin{equation}\label{eq:switching_con2}
	F_2(\p_2) \cdot \nabla H = -2 d_1 \eta X' - \eta^2 \, X - \frac{7}{2}\mu 
	\left[2d_2 Y' + Y + \ve_2^2 g\right] =0,
\end{equation}
As before, we see that such point can only occur if $\Lambda_N=\Ord{\eta}$. A 
similar argument shows that
$Y'>0$ at any such point $\p_2$. The argument then follows as before. The 
second directional derivative at the switching point will thus be negative, 
once again confirming consistency of Filippov's theory.

The key question for final consideration is whether the ball can lift off 
(that is, reach the termination point of bounce) while rolling. For that to be 
the case, we would need $\Lambda_N=0$ as the 
lift off condition, but thus also $\Lambda_T=0$, which can be thought of $F_s$ 
being tangent to 
$F_1$ and $F_2$ simultaneously. Considering the form of the function 
$\alpha(\p)$ we would have that precisely at lift off we must simultaneously 
satisfy the two conditions 
\begin{equation}
	4d_1 \eta X' + \eta^2 X=0 \quad \mbox{and} \quad 
	2d_2Y'+Y+\ve_2^2 g =0,  
	\label{eq:2conds}
\end{equation}
which imposes an additional constraint on the initial-value problem. Thus there 
will be at most a codimension-one surface in the
space of initial conditions for which both conditions can be 
simultaneously satisfied and lift off with rolling would occur with 
probability zero. However, to make this statement more precise requires a 
proper analysis of the singularity where $F_s$ is simultaneously tangent to 
$F_1$ and $F_2$, because it is possible that the point satisfying 
\eqref{eq:2conds} could be an attractor for sliding initial conditions. Such an 
analysis forms the subject of the next section. 

\begin{figure}[!htb]
	\centering
	\begin{subfigure}{0.48\linewidth}
		\centering
		\begin{overpic}[width = 0.8\linewidth]{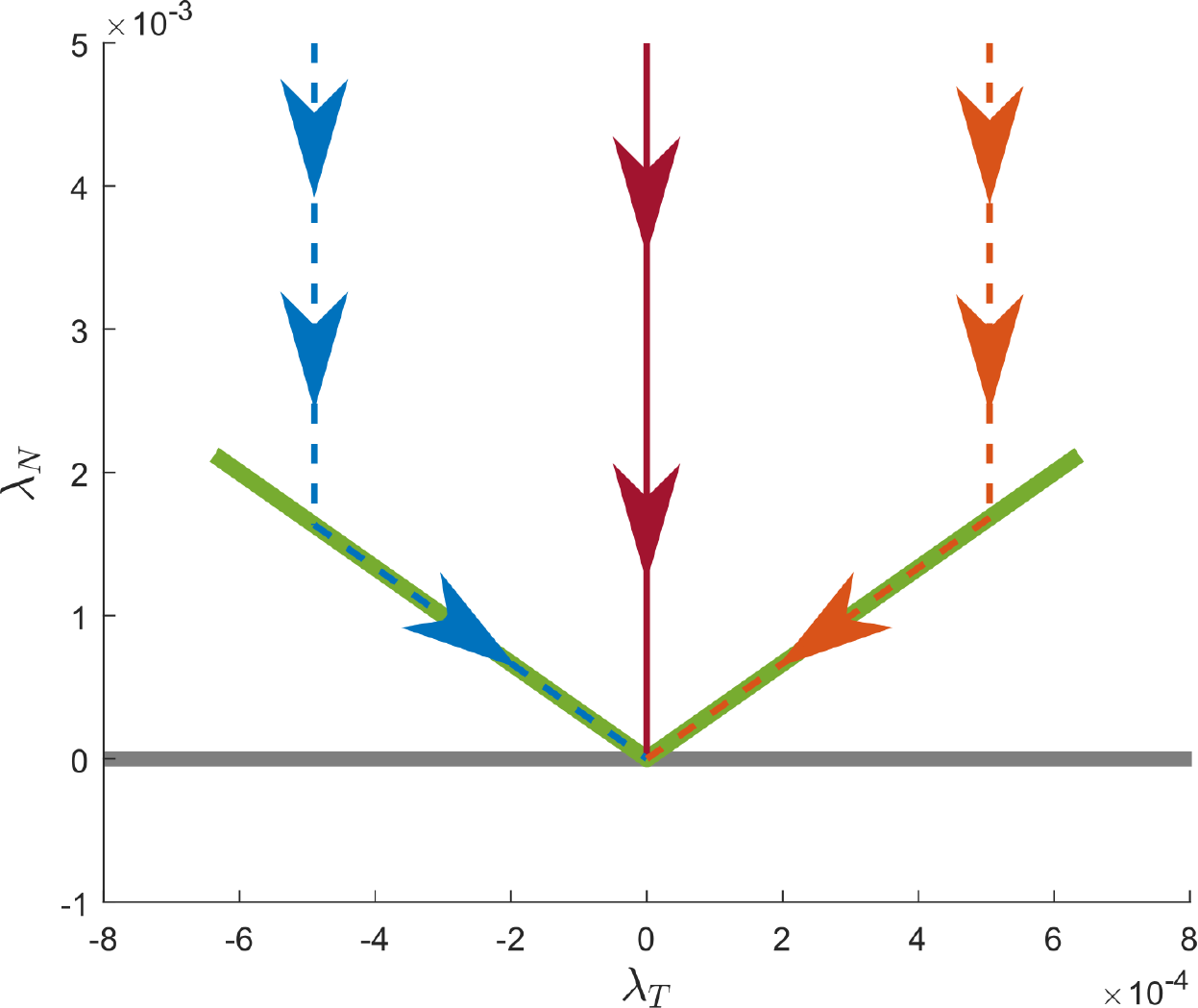}
			\put(42,72){(A)}
			\put(14,72){(B)}
			\put(71,72){(C)}
			\put(76,25){\small $\lambda_T = \mu \lambda_N$}
			\put(12,25){\small $\lambda_T = -\mu \lambda_N$}
			\put(15,13.5){Lift off}
		\end{overpic}
		\caption{ }
	\end{subfigure}
	\begin{subfigure}{0.48\linewidth}
		\centering
		\begin{overpic}[width = 0.8\linewidth]{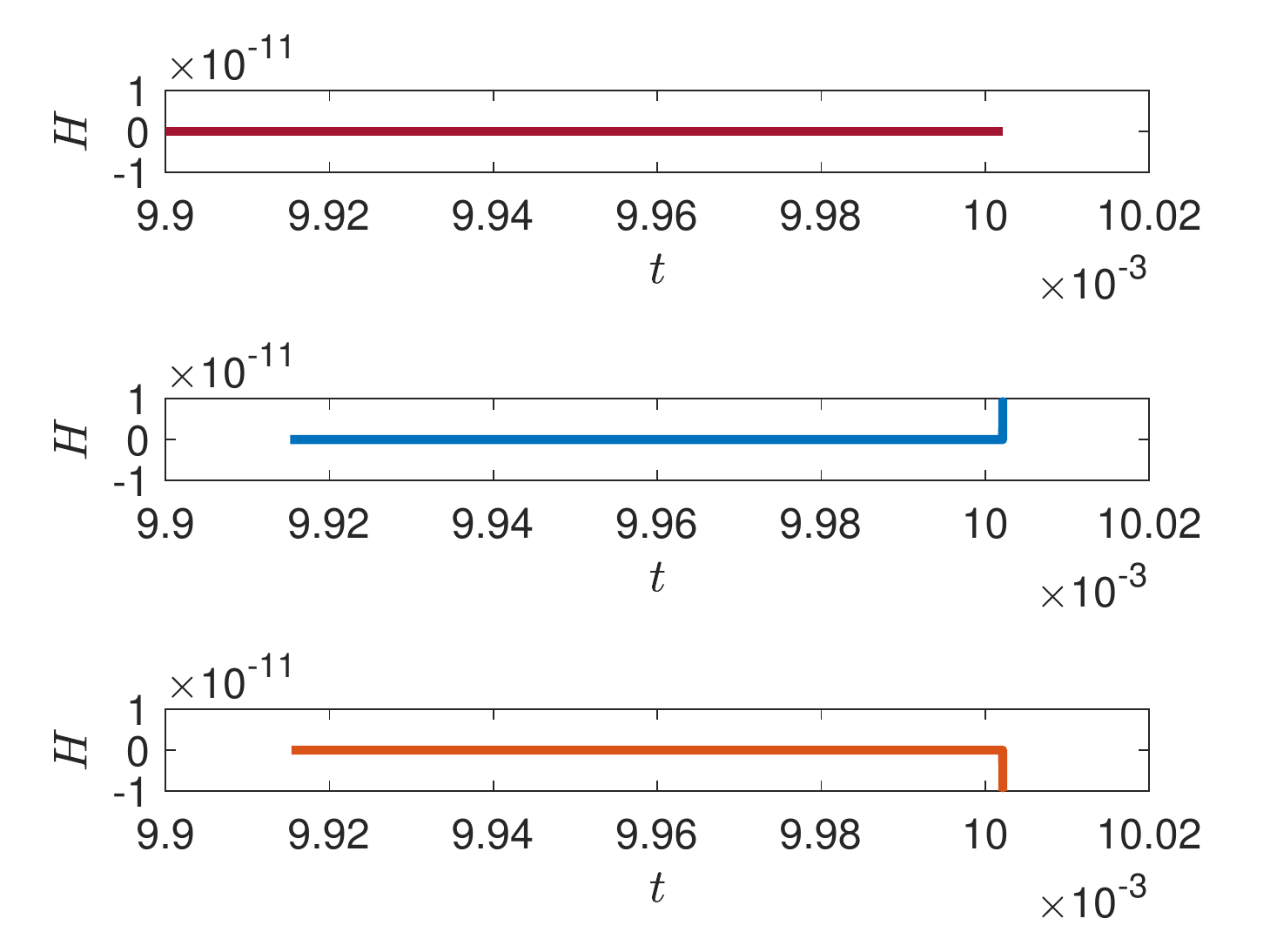}
			\put(11,81){(A)}
			\put(11,47){(B)}
			\put(11,15){(C)}
		\end{overpic}
		\caption{ }
	\end{subfigure}
	\caption{Numerical results for the system
		\eqref{eq:kv_rescaled} for $\mu=0.3$, with three different initial 
		conditions;  (a) projected onto the {\em friction cone} 
		$(\lambda_T,\lambda_N)$ and as slip velocity $H$ as a function of time. 
		See text for details.}
	\label{fig:kv_pert}
\end{figure}
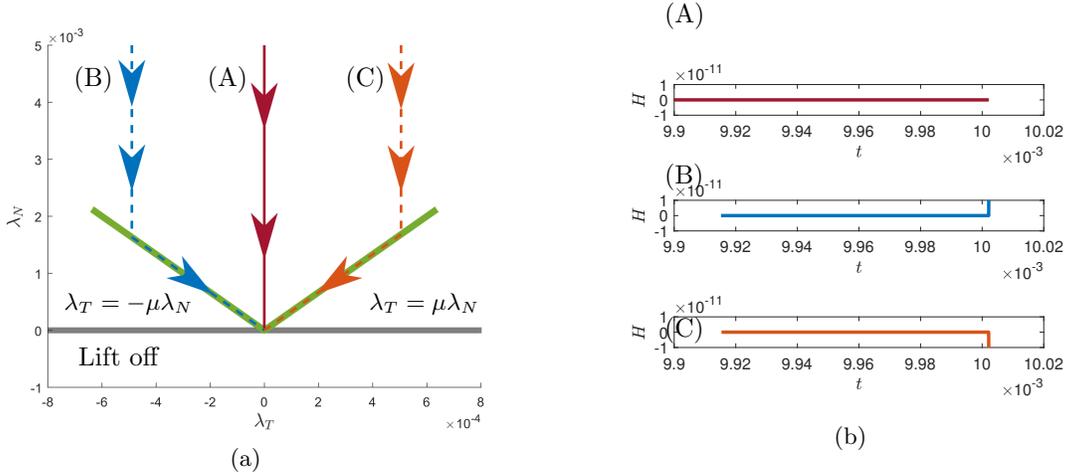

Briefly though, Figure \ref{fig:kv_pert} shows a numerical verification that 
lift-off with roll seems to not be an attractor in this case. Here trajectory 
(A) with initial conditions $\p_0 = [0.3543,\, -0.1603128, \,-0.1608, \, 
3.4739,\, 0.1603128]^{\intercal}$ is the trajectory which leaves the bounce 
within the friction cone, thus rolling. Initial conditions of trajectory (A) is 
perturbed slightly and gives raise to trajectories (B) and (C), which lift off 
on the boundary of the slipping cone, therefore slipping. (b) Tangential 
velocities of each of the trajectories. Note how trajectories (B) and (C) enter 
slipping motion (non-zero tangential velocity $H$) shortly before lift-off.

Thus we have found precisely the opposite conclusion from the rigid impact 
model in Section
\ref{sec:rigid_bounce}, namely that an introduction of even minimal, compliance 
in the normal 
direction causes the ball to always lift off slipping. That is, despite a 
transition to roll
during the bounce, with probability zero does does the ball lift off rolling.
This conclusion appears to be consistent with the experimental data, but it 
would be useful to check whether it is an artefact of the specific linear 
Kelvin-Voigt model we have chosen. 


\section{A generalised elasto-plastic model } \label{sec:analysis}

To test the generality of our findings to refinements in the surface model, we 
now consider a generalisation of the model form the previous section by posing
equations of motion of the form  
\begin{equation}\label{eq:generic_eq}
	\begin{aligned}
		\ddot{x} +  u(x,\dot{x}, y, \dot{y}, \omega) \, \dot{x} + z(x,\dot{x},
		y, \dot{y}, \omega)\,x & = \lambda_T,\\
		\ddot{y} +  d(y, \dot{y}) \, \dot{y} + k( y, \dot{y})\, y & = -g ,\\
		\dot{\omega} & = \frac{5}{2} \lambda_T,
	\end{aligned}
\end{equation}
where $u,\, z,\, d,\, k$ are general, possibly nonlinear, functions, $g$ is the 
acceleration due to gravity and all other variables have the same meaning as 
before. Note this model allows for nonlinear, and velocity dependent stiffness 
and damping coefficients, as well dependence of the tangential mechanics on 
depth, spin and normal  velocity. We will however impose some asymptotic 
constraints on the functions $u,\, z,\, d$ and $k$ in the course of our theory 
(see equations 
\eqref{eq:functions_constraints}-\eqref{eq:functions_constraints_horizontal_scaling}
below).

Under our assumption of a point of contact, we assume that the frictional 
contact
between the ball and the surface can be modelled 
according to the Coulomb friction law, that is 
$$
\left|\lambda_T\right| = \mu \lambda_N = \mu \left( -d(y,\dot{y}) \, \dot{y} 
-k(y,\dot{y})\, y -g\, \right)
$$
during slipping and 
$$
\left|\lambda_T\right| < \mu \lambda_N = \mu \left( -d(y,\dot{y}) \, \dot{y} 
-k(y,\dot{y})\, y -g\, \right)
$$
when the ball is rolling (sticking), where $\mu$ is the coefficient of 
friction. 

The ball will remain in contact with the surface until normal forces due to the 
mechanism are balanced by those due to gravity, that is when 
\begin{equation}\label{eq:lift_off_cond}
	\lambda_N = -d(y,\dot{y}) \, \dot{y} -k(y,\dot{y})\, y -g = 0
\end{equation}
We will refer to \eqref{eq:lift_off_cond} as the \textit{lift-off condition}; 
the
horizontal and vertical velocities at this instance will 
be denoted by $\dot{x}_F$ and $\dot{y}_F$ respectively, and the lift-off spin 
denoted
by $\omega_F$.

\begin{figure}[!htb]
	\centering
	\begin{overpic}[width = 0.4\linewidth]{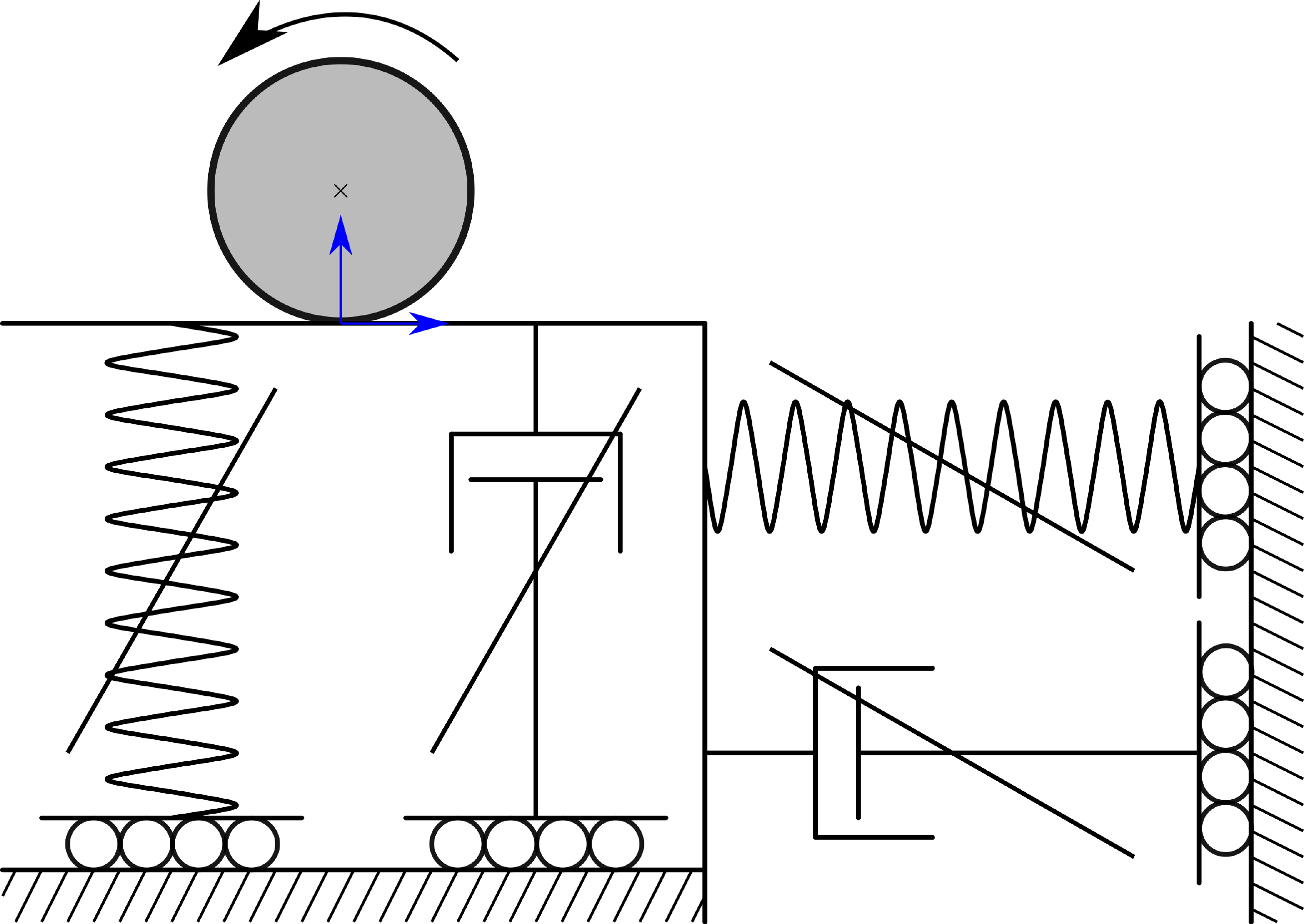}
		\put(22,58){\tiny $(x, y)$}
		\put(30,71){\tiny $\omega$}
		\put(42,20){\tiny $d(y,\dot{y})$}
		\put(-6,20){\tiny $k(y,\dot{y})$}
		\put(65,1){\tiny $u(x, \dot{x}, y,\dot{y}, \omega)$}
		\put(65,45){\tiny $z(x, \dot{x}, y,\dot{y}, \omega)$}
	\end{overpic}
	\caption{A generalised elasto-plastic point contact model.}
	\label{sfig:setting_generic}
\end{figure}
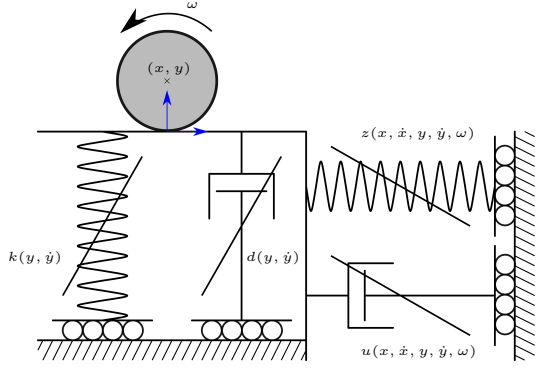

One can think of the form of \eqref{eq:generic_eq} as a Kelvin-Voigt-like 
model, where springs and
dampers, connected in parallel, dictating the behaviour of the sphere
during the bounce in the vertical and horizontal direction independently. The 
major 
difference, between the classical Kelvin-Voigt model and the one
used here is that the springs and dampers are allowed to be nonlinear and to 
feature coupling between normal and tangential degrees of freedom. 
An illustration of this setup can be seen in Figure \ref{sfig:setting_generic}.

We impose certain conditions on our functions $u,z,d,k$ suggested by physical
intuition. Firstly, we require the stiffness and damping in both vertical and
horizontal directions to increase as the ball moves further into the bouncing
surface, that is those values should increase as the vertical position 
$y$ decreases:
\begin{equation}\label{eq:functions_constraints}
	\begin{aligned}
		u_y (x,\dot{x}, y, \dot{y}, \omega)& < 0,\\
		d_y(y, \dot{y})  & < 0.
	\end{aligned}
\end{equation}
The motivation for this, is that in reality, point contact is an approximation 
to a Hertzian-like contact
in which the true area of contact increases with depth $-y$.
In a similar fashion, we require that the normal force acting on the body 
increases with depth
\begin{equation}\label{eq:force_constraints}
	\frac{\partial}{\partial y}\lambda_N = \frac{\partial}{\partial 
		y}\left(-d(y,\dot{y} ) \dot{y} - k(y,\dot{y}) y - g\right) <0. 
\end{equation}
Given that the ball is landing with $\dot{x}_0>0$, in the horizontal direction 
we will require that the stiffness increases when the ball is moving to the 
right, and decreases when the ball reverses (moves to the left):
\begin{equation}\label{eq:functions_constraints_horizontal}
	\sign{u_{\dot{x}} (x,\dot{x}, y, \dot{y}, \omega)} < 0 ) = \sign{\dot{x}}. 
\end{equation}

Our analysis is motivated by the case of a rapid ball bounce, which we suppose 
to take $\Ord{\ve}$ amount of time, where $\ve \ll 1$. We thus introduce a 
re-scaling of time
\begin{equation}
	\tau = \ve^{-1} t 
\end{equation}
and with it, we denote state variables in the new variable as 
$\ns{x}(\tau) = x(t),\, \ns{y}(\tau) = y(t),\, \nsomega(\tau) = 
\omega(t)$.
To ensure the appropriate balance of terms in \eqref{eq:generic_eq} we now 
require
\begin{equation}\label{eq:functions_constraints_horizontal_scaling}
	\begin{aligned}
		k(y, \dot{y}) & = \ve^{-1} \gk (\ns{y}, \ns{y}')\\
		d(y, \dot{y}) & = \ve^{-2} \gd (\ns{y}, \ns{y}'),
	\end{aligned}
\end{equation}
where $' = \frac{\dd}{\dd \tau}$ denotes the derivative with respect to the 
new time variable $\tau$, and $\gk (\ns{y}, \ns{y}'), \gd (\ns{y}, \ns{y}') 
\sim \Ord{1}$ are general functions.

At the same time, we anticipate a much greater stiffness  in the horizontal 
direction than in the vertical section. One can think of this requirement as 
avoiding the case of
a ``jelly-like'' surface, where the 
symmetry in vertical and horizontal direction would causes an outbound 
trajectory, in the absence of spin,
to be the time-reverse of the incoming trajectory, irrespective of the angle of 
incidence with the surface.
We will thus introduce a distinguished scaling so that the damping and 
stiffness function in the horizontal directions to be
\begin{equation}
	\begin{aligned}
		\gu (\ns{x}, \ns{x}', \ns{y}, \ns{y}', \nsomega) & = u(x, \dot{x}, 
		y, \dot{y}, \omega) \sim \Ord{1},\\
		\gz (\ns{x}, \ns{x}', \ns{y}, \ns{y}', \nsomega) & = z(x, \dot{x}, 
		y, \dot{y}, \omega) \sim \Ord{1}.
	\end{aligned}
\end{equation}

Our final requirement on the rescaled function $\gk$, $\gz$, $\gd$ and $\gu$ is 
that their (partial) 
derivatives with respect to the state variables be at most of $\mathcal{O}(1)$. 
That is 
we do not expect abrupt changes in the material constants with displacement or 
velocity.
That is, taking $q$ to be any 
of the state variables $\ns{x}, \ns{x}', \ns{y}, \ns{y}', \nsomega$ we have 
that
\begin{equation} \label{eq:derivatives_constraint}
	\frac{\partial}{\partial q} \gu, \frac{\partial}{\partial q} \gz, 
	\frac{\partial}{\partial q} \gd, \frac{\partial}{\partial q} \gk \sim 
	\Ord{1}.
\end{equation}

Following the rescaling and the introduction of new functions, the resulting 
dynamical system for modelling the bounce becomes
\begin{equation}
	\begin{aligned} \label{eq:generic_rescaled}
		\ns{x}'' + \ve \gu(\ns{x}, \ns{x}', \ns{y}, \ns{y}' 
		\nsomega) \, \ns{x}' + \ve^2 \gz (\ns{x}, \ns{x}', \ns{y}, \ns{y}' 
		\nsomega) \, \ns{x} & = \Lambda_T,\\
		\ns{y}'' + \gd(\ns{y},\ns{y}')\, \ns{y}' + \gk(\ns{y}, \ns{y}') \, 
		\ns{y} & = -\ve^2 g,\\
		\nsomega' &= \frac{5}{2} \Lambda_T,	
	\end{aligned}
\end{equation}
where $\Lambda_T = \ve^2 \lambda_T $. 

The conditions introduced in \eqref{eq:functions_constraints}, 
\eqref{eq:force_constraints} and \eqref{eq:functions_constraints_horizontal} 
can now be written in the form
\begin{equation}\label{eq:functions_constraints_new_scaling}
	\begin{aligned}
		\gu_{\ns{y}} (\ns{x},\ns{x}', \ns{y}, \ns{y}', \nsomega)& < 0,\\
		\gd_{\ns{y}} ( \ns{y}, \ns{y}')  & < 0,\\
		\sign{\gu_{\ns{x}'} (\ns{x},\ns{x}', \ns{y}, \ns{y}', \nsomega) } & = 
		\sign{\ns{x}'}
	\end{aligned}
\end{equation}
and 
\begin{equation}\label{eq:force_constraints_new_scaling}
	\frac{\partial}{\partial Y}\Lambda_N = \frac{\partial}{\partial 
		y}\left(-\gd(\ns{y},\ns{y}')\, \ns{y}' - \gk(\ns{y}, \ns{y}') \, 
	\ns{y} -\ve^2 g \right) <0. 
\end{equation}

The Coulomb friction law once again defines our dynamics for the slipping 
motion, 
where
\begin{equation}
	\Lambda_T = - \sign{H(\p)} \Lambda_N,
\end{equation}
where
\begin{equation}
	H(\p) = \ns{X}' + \nsomega \quad \mbox{and} \quad\Lambda_N = 
	-\left( \gd \ns{y}' + \gk Y + \ve^2 g\right),
\end{equation}
and $\p$ is once again the five-dimensional vector of dynamical coordinates 
defined by \eqref{eq:vec_p}. Compared to the simple analysis in Section 
\ref{sec:two_cases} additional consideration is needed for 
the rolling motion, when $H(\p) = \ns{X}' + 
\nsomega=0.$ Precise study of the onset and
loss of rolling motion and their implication for the overall dynamics will be 
our main interest in what follows. 

\subsection{The two-fold singularity}

Let us consider again the possibility of  lift off whilst rolling. Recall that 
at lift off $\Lambda_N = 0$ 
and as such the definition of friction dictates that we must also have 
$\Lambda_T=0$. At the same time, to agree with  the formulation as a Filippov 
System, the lift off point must coincide with the lower dimensional spaces 
where $\ga(\p) = 0$ and $\ga(\p) = 1$. Such point is denoted as $\s$  
in Figure \ref{fig:geometries}.

The point where the two switching surfaces ($\ga = 0$ and $\ga =1$) coincide is 
known in the literature as the two-fold singularity \cite{jeffrey_update}. It 
is now of our interest to determine whether such point is attracting or not. 
Indeed, an attracting two-fold singularity (which we will from now on denote as 
point $\s$), will mean that rolling trajectories will tend towards it, 
therefore reaching lift off with $H(\p)=0$. Alternatively, if $\s$ is not 
attracting, there will only be a finite, number of lower-dimensional 
trajectories tending to $\s$.

%
%
A full unfolding of two-fold singularity can be found in \cite{jeffrey_hidden}
with the formal results following from Theorems 6.1 \& 6.2 in the 
aforementioned book. We look at the case of our interest only in this paper.

Let us consider an $n$-dimensional discontinuous system
\begin{equation*}\label{eq:generic_filippov}
	\p' = \begin{cases}
		F_1(\p) & \quad \mbox{if } \quad H(\p) > 0,\\
		F_2(\p) & \quad \mbox{if } \quad H(\p) < 0,
	\end{cases}
\end{equation*} 
where $H(\p) = 0$ is the discontinuous surface, along which $\p' = F_s(\p) = (1 
- \ga(\p))\, F_1(\p) + \ga(\p) \, F_2(\p)$, with $\ga$ defined by Equation 
\eqref{eq:alpha} as before. 

Suppose the system posses two switching surfaces $\ga(\p) = 0$ (equivalent with 
$(F_1 \cdot \nabla) H(\p)$) and $\ga(\p) = 1$ (equivalent with $(F_2 \cdot 
\nabla)H(\p)$) which intersect at a point $\s$.

We define new coordinates:
\begin{equation}
	\begin{aligned}
		z_1 & = \gt (\p) H(\p),\\
		z_2 & = - \gp(\p) F_1 \cdot \nabla \left( \gt(\p) H(\p) \right),\\
		z_3 & = \frac{F_2 \cdot \nabla \left( \gt(\p) H(\p) \right)}{\gp (\p)},
	\end{aligned}
\end{equation}
where
\begin{equation}
	\gt(\p) = \left|(F_1\cdot \nabla)^2 H(\p) \, (F_2\cdot \nabla)^2 H(\p) 
	\right|^{-1/2}, \qquad \gp (\p) = \left|\frac{(F_2\cdot \nabla)^2 H(\p) 
	}{(F_1\cdot \nabla)^2 H(\p)} \right|^{1/4}.
\end{equation}

We note that $z_1=0$ is precisely the discontinuity surface of 
\eqref{eq:generic_filippov} with $z_2=0$ being equivalent to the switching 
surface $\ga=0$ and $z_3$ being equivalent to $\ga =1$. Thus, the two-fold 
singularity is precisely the $n-3$ dimensional manifold $\{(\z=z_1, \dots, z_n) 
\,:\, z_1 = z_2 = z_3 = 0\}$, with $z_{4},\dots, z_n$ defined to be any $n-3$ 
dimensional coordinates orthogonal to $z_1, z_2$ and $z_3$.

In deriving the normal form of the dynamics we apply a time rescaling $t 
\mapsto \gp t$ when $z_2<$ and $t \mapsto t / \gp$ when $z_2>0$. Since $\gp>0$ 
the phase portrait remains unchanged under this transformation. 

With details and series expansions presented in the original source, one finds 
that in the neighbourhood of the two-fold singularity, under the 
transformations outlined before, the flow becomes
\begin{equation}
	\begin{aligned}
		\relax
		[\dot{z}_1, \dot{z}_2, \dot{z}_3] = &
		\begin{Bmatrix}
			\left[- z_2,\, -\sigma_1,  
			\nu_{1} \right] & 
			\quad \mbox{if} \quad z_1 >0\\
			\left[z_3, \, \nu_{2}, \, \sigma_2 
			\right] & 
			\quad \mbox{if} \quad z_1 <0\\
		\end{Bmatrix} + \left[\Ord{|\z|^2},\Ord{|\z|},\Ord{|\z|}\right],\\
		\dot{z}_{4,5,\dots,n} & = \Ord{|\z|},
	\end{aligned}
\end{equation}
where 
\begin{equation} \label{eq:sigma_def}
	\sigma_i = \sign{\left(F_i \cdot \nabla\right)^2 H(\s)}
\end{equation}
and
\begin{equation} \label{eq:nus}
	\gv^{1} = \left.\frac{\left(F_1 \cdot \nabla\right) \left(F_2 \cdot 
		\nabla\right) H }{\sqrt{\left|\left(F_1 \cdot \nabla\right)^2 H \cdot 
			\left(F_2 \cdot \nabla\right)^2 H\right|}} \right|_{\z=\s}, \qquad
	\gv^{2} = \left.-\frac{\left(F_2 \cdot \nabla\right) \left(F_1 \cdot 
		\nabla\right) H }{\sqrt{\left|\left(F_1 \cdot \nabla\right)^2 H \cdot 
			\left(F_2 \cdot \nabla\right)^2 H\right|}} \right|_{\z=\s}.
\end{equation}
As discussed in \cite{jeffrey_update} the quantities $\gv_{1,2}$ characterise 
the local curvature of the flow, with $\gv_1 \gv_2$ quantifying the jump in the 
vector field $F_{1,2}$ and the singularity.

With the new variables, along the discontinuity surface $z_1 = 0$ we define (as 
before)

\begin{equation}\label{eq:normal_sliding}
	[\dot{z}_1, \dot{z}_2, \dot{z}_3] = \left(1-\ga\right) \left[- z_2,\, 
	-\sigma_1, \nu_{1} \right] + \ga \left[_3, \, \nu_{2}, \, \sigma_2 
	\right] + [\Ord{|\z|^2},\Ord{|\z|},\Ord{|\z|}],
\end{equation}
where $\ga$ is solved so that $\dot{z}_1=z_1 = 0 $. To the leading order the 
dynamics in \eqref{eq:normal_sliding} are independent of $z_{4,\dots,n}$. 
Therefore, the dynamics in $z_2$ and $z_3$ can be separated and, with $\ga$ 
evaluated explicitly, we have the two-dimensional system
\begin{equation} \label{eq:normal_red}
	\begin{bmatrix}
		\dot{z}_2\\ \dot{z}_3
	\end{bmatrix} = \frac{1}{z_2+z_3}
	\begin{bmatrix}
		\gv_2 & - \sigma_1\\ \sigma_2 & \gv_1
	\end{bmatrix}
	\begin{bmatrix}
		z_2 \\ z_3
	\end{bmatrix}.
\end{equation}

The two switching surfaces (or \textit{folds}) $z_2$ and $z_3$ will divide the 
vector field along the discontinuity $z_1=0$ into four regions
\begin{equation}
	\begin{aligned}
		&\Da = \left\{ \z : z_1 = 0; z_2, z_3>0 \right\}, \\
		&\Dr = \left\{ \z : z_1 = 0; z_2, z_3<0 \right\}, \\
		&\Dc{1} = \left\{ \z : z_1 = 0; z_2 < 0 < z_3 \right\}, \\
		&\Dc{2} = \left\{ \z : z_1 = 0; z_3 < 0 < z_2 \right\}.
	\end{aligned}
\end{equation}
With such definition, the vector field outside of the discontinuity $z_1=0$ is 
directed into $\Da$, away from $\Dr$ and crosses $\Dc{1}$ in the direction of 
increasing $z_1$, and crosses $\Dc{2}$ in the direction of decreasing $z_1$.

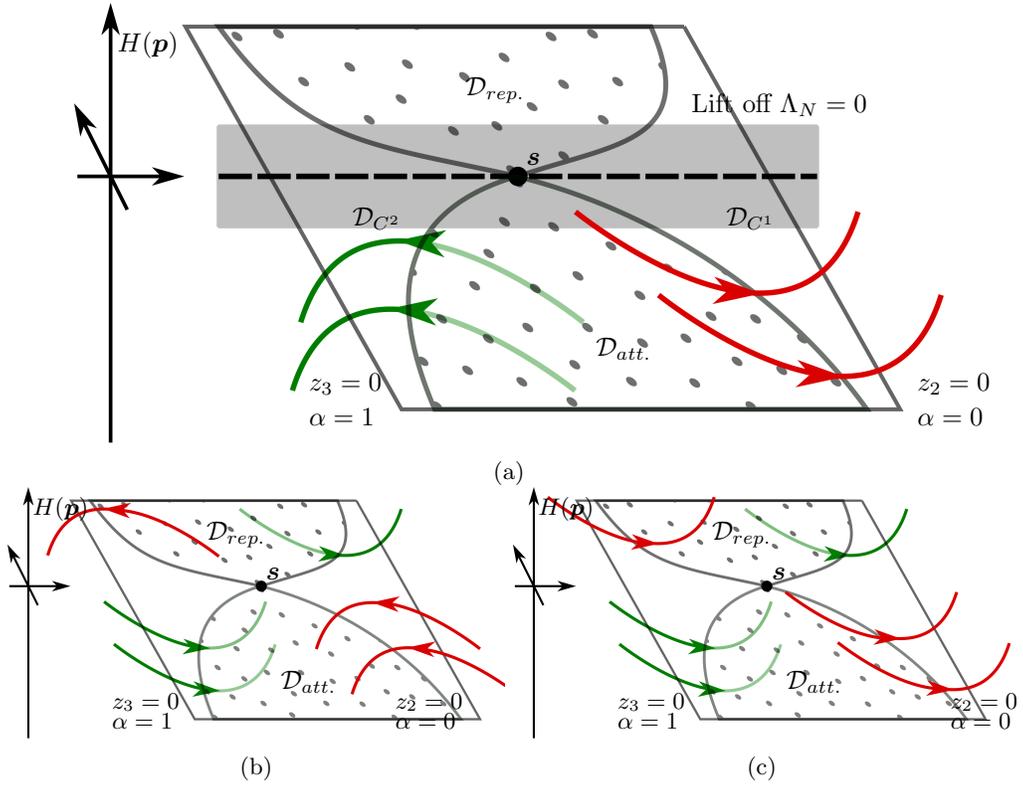
\begin{figure}
	\centering
	\begin{subfigure}{0.7\linewidth}
		\begin{overpic}[width = \linewidth]{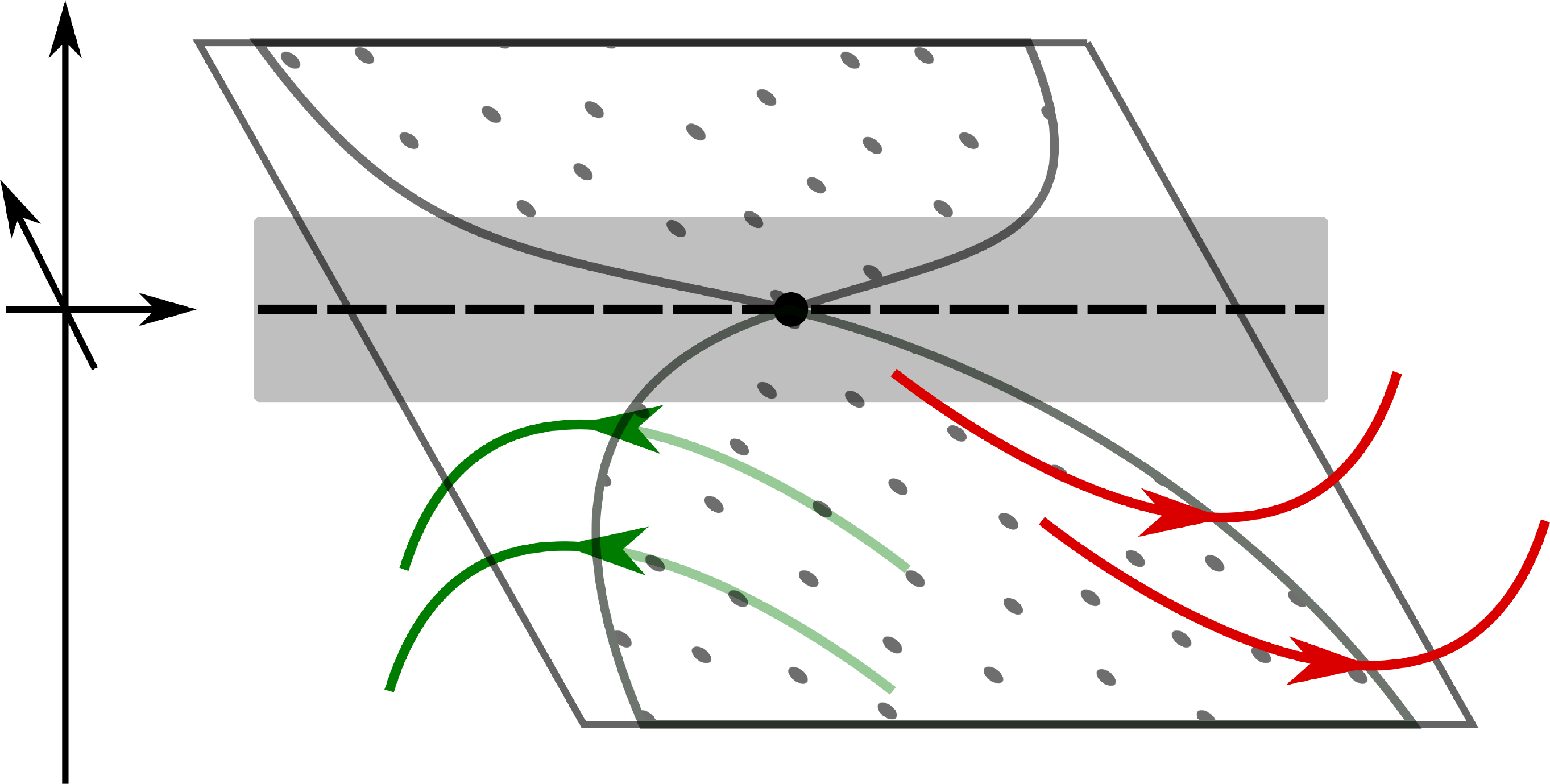}
			\put(5,45){$H(\p)$}
			\put(71,38){Lift off $\Lambda_N=0$}
			\put(52,32){$\s$}
			\put(27,6){$z_3=0$}
			\put(27,2){$ \alpha = 1 $}
			\put(97,6){$z_2=0$}
			\put(97,2){$ \alpha = 0 $}
			\put(60,10){$\Da$}
			\put(45,40){$\Dr$}
			\put(32,25){$\Dc{2}$}
			\put(75,25){$\Dc{1}$}
		\end{overpic}
		\caption{}
		\label{sfig:visible-2-fold}
	\end{subfigure}
	\begin{subfigure}{0.4\linewidth}
		\begin{overpic}[width = \linewidth]{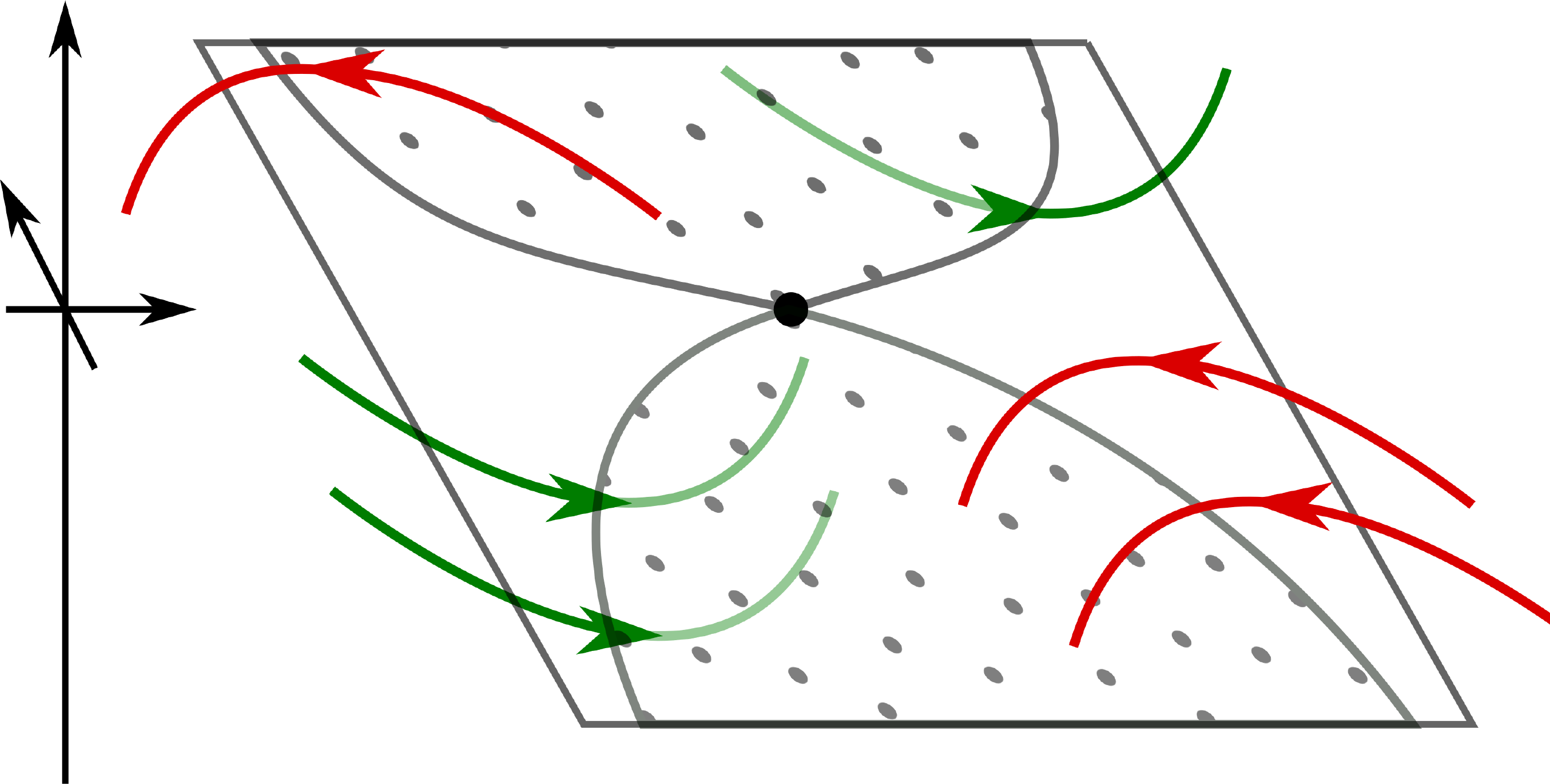}
			\put(5,45){\small $H(\p)$}
			\put(52,32){$\s$}
			\put(21,6){\small $z_3=0$}
			\put(21,2){\small $ \alpha = 1 $}
			\put(78,6){\small $z_2=0$}
			\put(78,2){\small $ \alpha = 0 $}
			\put(55,10){$\Da$}
			\put(40,40){$\Dr$}
		\end{overpic}
		\caption{}
		\label{sfig:invisible-2-fold}
	\end{subfigure}
	\begin{subfigure}{0.4\linewidth}
		\begin{overpic}[width = \linewidth]{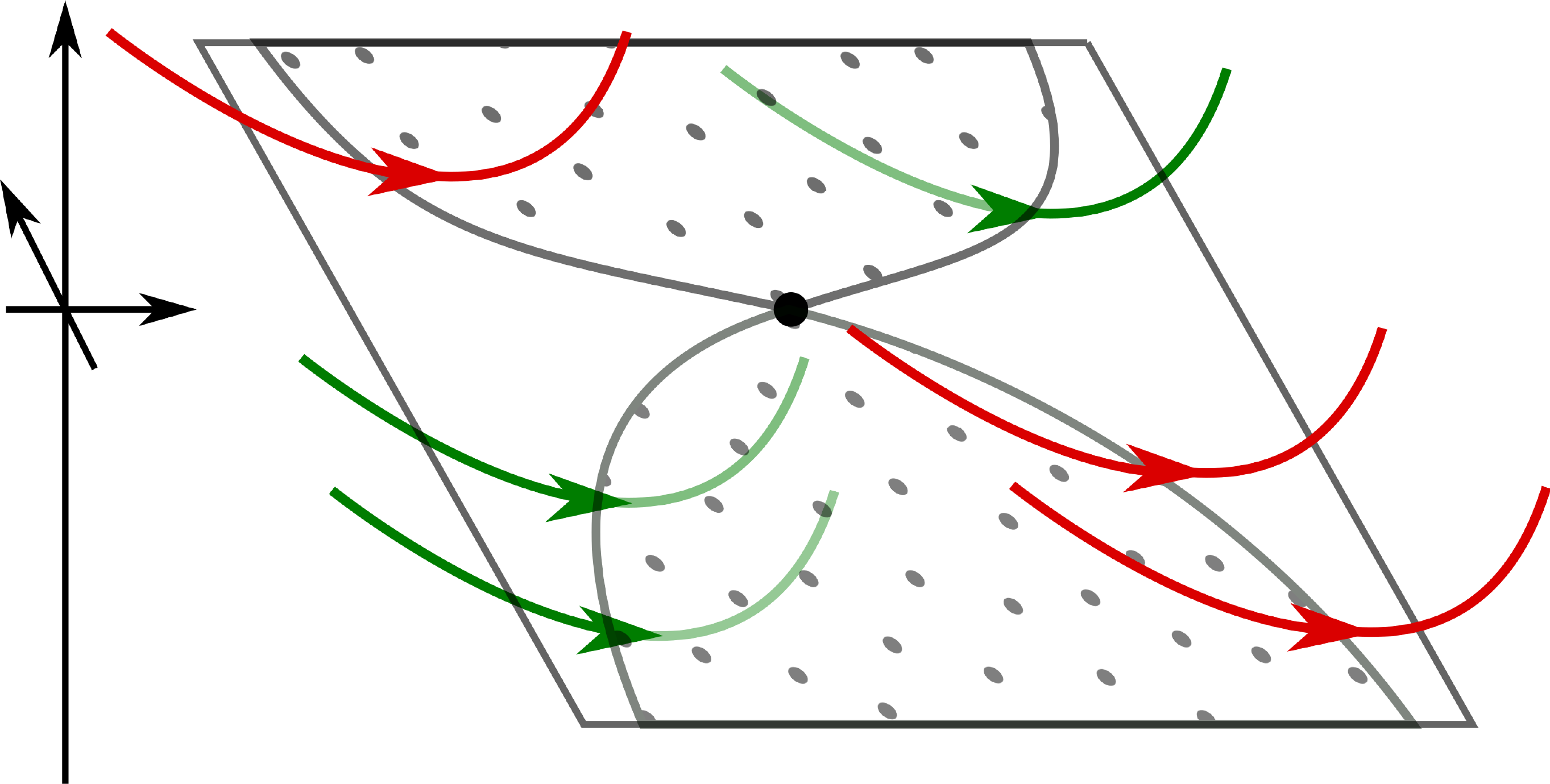}
			\put(5,45){\small $H(\p)$}
			\put(52,32){$\s$}
			\put(21,6){\small $z_3=0$}
			\put(21,2){\small $ \alpha = 1 $}
			\put(88,6){\small $z_2=0$}
			\put(88,2){\small $ \alpha = 0 $}
			\put(55,10){$\Da$}
			\put(40,40){$\Dr$}
		\end{overpic}
		\caption{}
		\label{sfig:in-visible-2-fold}
	\end{subfigure}	
	\caption{Classification of two-fold singularities, depending on the flow 
		outside of the discontinuity surface. (a) Visible two-fold (interpreted 
		for the 
		case of a ball bounce, with the lift-off surface indicated). Here, red 
		trajectories correspond to flow with $H>0$ and green trajectories with 
		$H<0$. (b,c) Adaptions to the above figure to deal with the (b) 
		invisible 
		two-fold and (c) the visible-invisible two-fold. Figure redrawn based 
		on images in 
		\cite{jeffrey_update}.}
	\label{fig:2-fold}
\end{figure}

According to \cite{jeffrey_hidden}, three types of two-fold singularity can be 
distinguished based on the nature of the flow outside of the discontinuity 
surface, specifically the nature of the {\em grazing} tangencies at the sliding 
region; see Fig.\ref{fig:2-fold}. Each of these cases are analysed in detail, 
here we analyse only the case that is relevant to the ball-bounce problem, 
namely the visible case. 

\subsection{Visible two-fold singularity}

When $\sigma_1>0 > \sigma_2$ then $\s$ is the intersection of visible folds, 
and we call such point a \textit{visible two-fold}. Trajectories away from the 
discontinuity are presented in Figure \ref{sfig:visible-2-fold}.
In this case \eqref{eq:normal_red} becomes
\begin{equation} \label{eq:normal_red_visible}
	\begin{bmatrix}
		\dot{z}_2\\ \dot{z}_3
	\end{bmatrix} = \frac{1}{z_2+z_3}
	\begin{bmatrix}
		\gv_2 & - 1\\ -1 & \gv_1
	\end{bmatrix}
	\begin{bmatrix}
		z_2 \\ z_3
	\end{bmatrix}.
\end{equation}

Understanding the nature of the fixed point $(0,0)$ is a straightforward, 
albeit laborious, exercise in algebra. Detailed derivation can be found in 
\cite[Chapter~13.4]{jeffrey_hidden}, we give the simple summary here.

\begin{figure}
	\centering
	\begin{subfigure}{0.3\linewidth}
		\begin{overpic}[width = \linewidth]{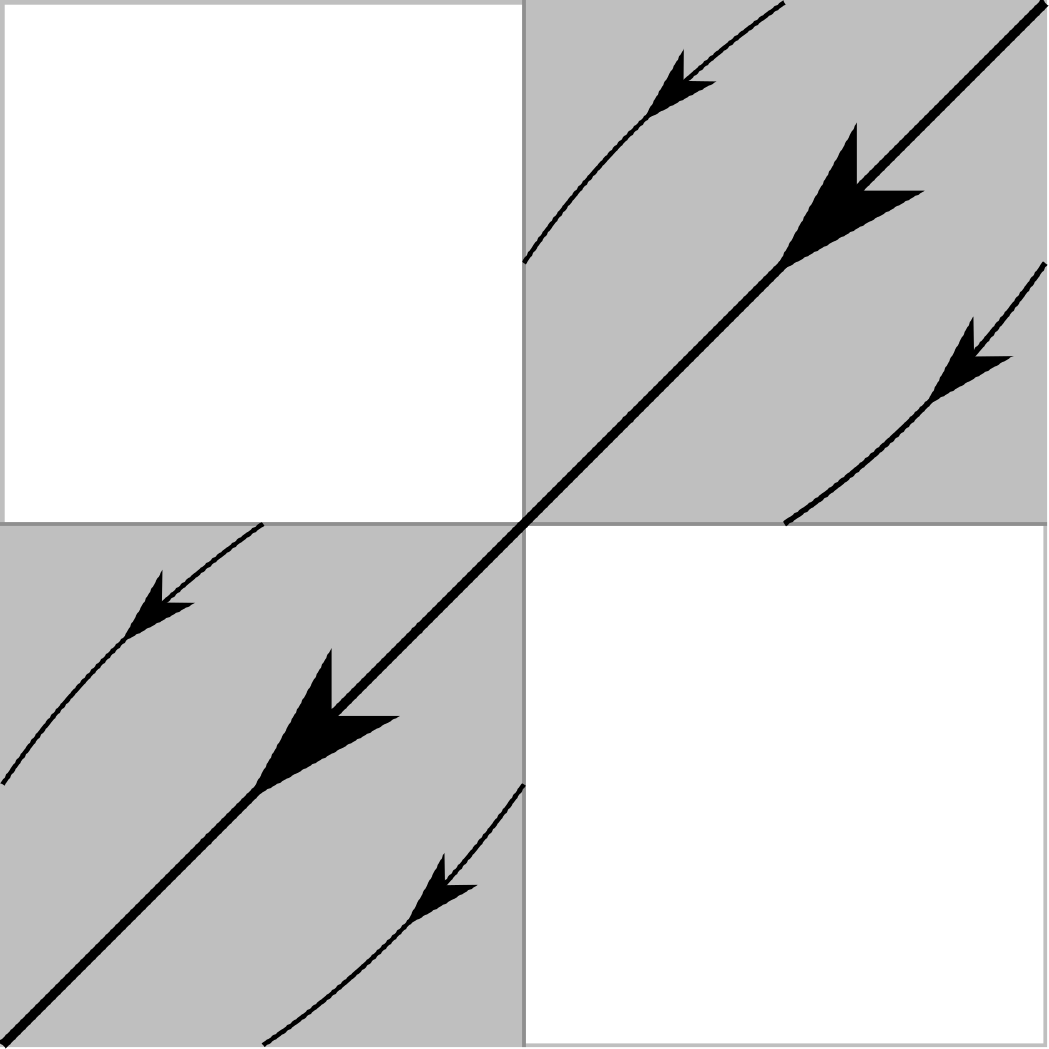}
			\put(82,72){$\Da$}
			\put(5,20){$\Dr$}
			\put(5,72){$\Dc{1}$}
			\put(82,20){$\Dc{2}$}
			\put(45,51){$\s$}
		\end{overpic}
		\caption{}
		\label{sfig:visible-fp1}
	\end{subfigure}
	\hspace{2cm}
	\begin{subfigure}{0.3\linewidth}
		\begin{overpic}[width = \linewidth]{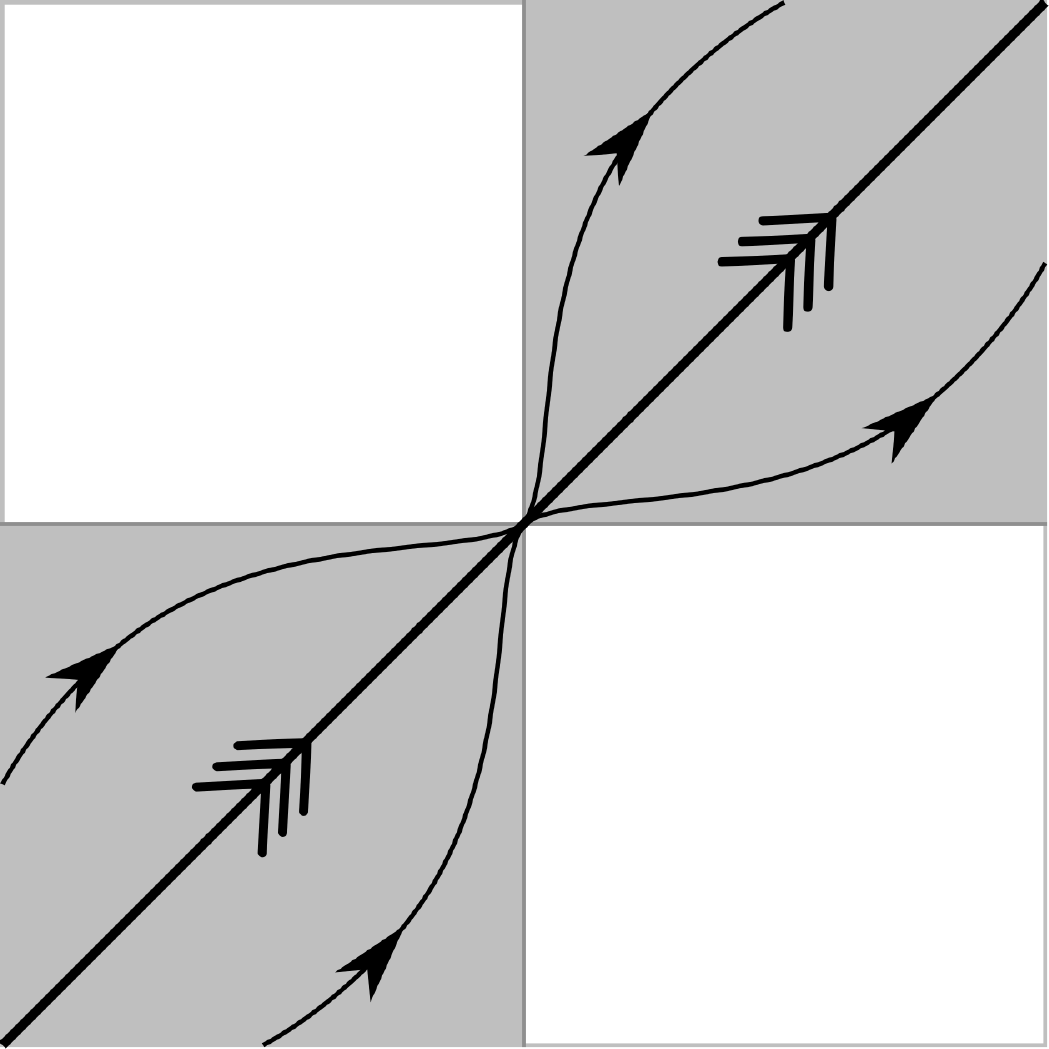}
			\put(82,72){$\Da$}
			\put(5,20){$\Dr$}
			\put(5,72){$\Dc{1}$}
			\put(82,20){$\Dc{2}$}
			\put(45,51){$\s$}
		\end{overpic}
		\caption{}
		\label{sfig:visible-fp2}
	\end{subfigure}
	\caption{Phase portraits for the visible two-fold singularity. (a) case 
		where $\gv_1 \gv_2<1$ or when $\gv_{1,2}<0$ and $\gv_1 \gv_2 > 1$. 
		Reproduced from \cite{jeffrey_update}.}
\end{figure}
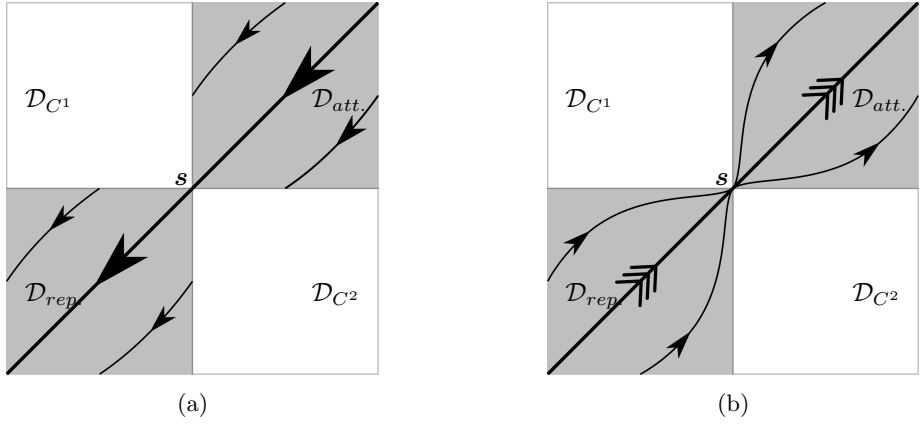

The linearised system \eqref{eq:normal_red_visible} will yield two eigenvalues 
$\gamma_+$ and $\gamma_-$.
For $\gv_{1,2}>0$ and $\gv_1 \gv_2>1$ we have that the two eigenvalues are such 
that $0< \gamma_- < \gamma_+$ and the weak eigendirection lies along 
eigenvector 
$\vec \gamma_-$ (corresponding to the eigenvalue $\gamma_-$). The flow will 
therefore generate a node, with the flow outwards from the singularity $\s$ 
and into the attracting region $\Da$. The phase portrait for that is shown in 
Figure \ref{sfig:visible-fp2}.

For $\gv_1 \, \gv_2 <1$ we have that $\gamma_-<0<\gamma_+$  and the inward 
eigendirection $\vec \gamma_-$ lies in the regions of discontinuity $\Da$ and 
$\Dr$. Only along the direction of $\vec \gamma_-$ we have a single trajectory 
reaching $\s$ from $\Da$ and all remaining trajectories will leave through the 
folds. The respective phase portrait is shown in Figure \ref{sfig:visible-fp1}.

When $\gv_{1,2}<0$ and $\gv_1 \, \gv_2>1$ we find that the eigenvalues are such 
that $\gamma_-<\gamma_+<0$ and again only a single trajectory along $\vec 
\gamma_-$ passes through $\Da$ to $\Dr$ through $\s$, and again we obtain the 
phase portrait from Figure \ref{sfig:visible-fp1}.

To proceed we therefore need to compute the normal form for the problem at hand
and determine the signs of
the quantities $\gv_1$ and  $\gv_2$.
Referring back to Figure \ref{fig:geometries} we see that reaching $\ga = 0$ 
(and thus entering $H(\p) >0$), whilst being tangential to the flow in $F_1$, 
the trajectory is also at its minimum with respect to $H(\p)$. A direct 
calculation shows that at $\s$ (and thus also in the neighbourhood of $\s$):
\begin{equation} \label{eq:second_derivativeF1}
	\begin{aligned}
		\left(\nabla\cdot F_1\right)^2 H(\s) &= \left.  \frac{7}{2} \mu \left[ 
		\gd_{\ns{y}} \, \ns{y}' + \gk_{\ns{y}} \, \ns{y} + \gk \right] \ns{y}' 
		+  \frac{7}{2} \mu \left[ 	\gd_{\ns{y}'} \, \ns{y}' + \gd 
		+ \gk_{\ns{y}'}\, \ns{y}  \right] \left(-\gd - \gk - \ve^2 g \right) 
		+ \Ord{\ve}\right|_{\p=\s}\\
		& = -\frac{7 \mu}{2} \frac{\partial}{\partial Y}\left(\Lambda_N\right) 
		Y' + \Ord{\ve},
	\end{aligned}
\end{equation}
where we have used that $\Lambda_N(\s) = 0$ together with the constraint on the 
size of the derivatives \eqref{eq:derivatives_constraint}. Noting that $\s$ can 
only be reached in the restitution phase (i.e. when $\ns{y}'>0$) we see that 
$\left(\nabla\cdot F_1\right)^2 H(\s)>0$, confirming thus that a trajectory 
leaving the rolling surface $H=0$ through the switching region $\ga=0$ will 
indeed enter the $F_1$ vector field as expected. A similar direct calculation 
shows that 
\begin{equation} \label{eq:second_derivativeF2}
	\left(\nabla\cdot F_2\right)^2 H(\s) = \frac{7 \mu}{2} 
	\frac{\partial}{\partial Y}\left(\Lambda_N\right) Y' + \Ord{\ve},
\end{equation}
which in turns confirms a consistent behaviour when leaving the rolling surface 
$H=0$ through the switching surface $\ga=1$. 

Recalling that the classification of the two-fold singularity depends on the 
quantities $\sigma_{1,2}$ (as defined by Equation \eqref{eq:sigma_def}) from 
Equations \eqref{eq:second_derivativeF1} and \eqref{eq:second_derivativeF2} we 
find that 
\begin{equation}
	\begin{aligned}
		\sigma_1 & = \sign{\left(\nabla\cdot F_1\right)^2 H(\s)} \\
		& = \sign{-\frac{7 \mu}{2} \frac{\partial}{\partial 
				Y}\left(\Lambda_N\right) Y' + \Ord{\ve}} = 1\\
		\sigma_2 & = \sign{\left(\nabla\cdot F_2\right)^2 H(\s)} \\
		& = \sign{\frac{7 \mu}{2} \frac{\partial}{\partial 
				Y}\left(\Lambda_N\right) Y' + \Ord{\ve}} = -1.
	\end{aligned}
\end{equation}
We thus have $\sigma_1>0>\sigma_2$, which indeed confirms from
\eqref{eq:generic_eq} that this is a 
is a visible two-fold singularity.

Let us now determine $\gv_{1,2}$. From \eqref{eq:nus} 
we have that
\begin{equation} \label{eq:nus_calc}
	\gv_1\,\gv_2 = 1 + \frac{20}{7} \, \frac{\frac{\partial}{\partial X'}(\gu) 
		X'}{\frac{\partial}{\partial Y'} (\Lambda_N ) Y'} g \ve^3 + 
	\Ord{\ve^4}.
\end{equation}
By our assumption \eqref{eq:force_constraints_new_scaling}
the denominator in \eqref{eq:nus_calc} is negative. Furthermore, the condition 
introduced in \eqref{eq:functions_constraints_new_scaling} 
dictates that $\frac{\partial}{\partial X'}(\gu) X'>0$, therefore yielding 
$\gv_1\,\gv_2<1$, and 
thus we obtain a phase portrait shown in Figure \ref{sfig:visible-fp1}.
This is the case where $\s$ is a non-attracting point, with 
only a single trajectory passing through it within the normal form. 

Interpreting our result with respect to the ball bounce, we note that the ball 
that enters rolling during the bounce (governed by the 
aforementioned  conditions) can only leave with a roll along a codimension-one 
manifold within the space of initial conditions. Thus, with probability 1, the 
ball must leave slipping.

\subsection{Rigid bounce revisited}

Let us briefly return to the case of the rigid bounce studied in
\ref{sec:rigid_bounce}. We concluded the studies in that section with
a result, where a ball that entered rolling during the bounce could no
longer enter the slip. One could think that this contradicts our main
result here --- but the case is more subtle than that.

The rigid bounce is a limiting case of the model studied in this paper, where 
the tangential compliance tends to zero;  $\ve \rightarrow 0$. This, in turn, 
means that the quantity $\gv_1\,\gv_2= 
1$, and hence the normal form of the visible two-fold singularity given by 
\eqref{eq:normal_red_visible} is based on a singular matrix. 
Classification is thus more convoluted than that for the general case, but has 
already been studied in detail in \cite[chapter~13]{jeffrey_hidden} and is  
known as the diabolo bifurcation.

Rather than studying the singularity for that case in detail, one can notice  
interesting geometries about that case. Since the lack of tangential compliance 
affects  the tangential force during the roll, the Filippov formulation should 
lead to 
an ``equal'' contribution from the two vector fields $F_1$ and $F_2$. Explicit  
calculation of the function $\ga(\p)$ as given by Equation \eqref{eq:alpha}  
shows $\ga=1/2$ throughout the roll, and thus $F_s(\p) = \left(F_1(\p) + F_2  
(\p)\right)/2$. 

The geometry on the friction cone is thus quite simple -- all trajectories 
that  enter the discontinuity (roll) region land precisely on the attracting  
trajectory that leads to a lift off with roll.  In other words -- that 
attracting trajectory is the only trajectory allowed in the rolling region for  
the case of the rigid bounce under Filippov's formulation.


\section{Conclusion} \label{sec:conclusion}

Inspired by empirical observations on the experimental data from 
\cite{experiments}, we used Filippov theory to analyse ball bounce within 
generalised point contact models that have both normal and tangential 
compliance. Under physically realistic scaling laws, we show that  a rigid 
spinning sphere will typically lift off with some non-zero relative tangential 
velocity. The case of rolling lift off forms the boundary case between forward 
and backward slipping cases, and therefore should not be observed in practice. 
This is contrary to a much more well studied and well understood problem of a 
rigid ball bouncing off a rigid surface. We have identified the underlying 
principle of the slip and roll transition, which has not been considered in 
detail before, and a careful understanding of that will be paramount to future 
modelling.

We are now in possession of what is thought to be the largest data set on the 
given problem, however fitting the models to the data remains to be an issue. 
Indeed, the focus of our analysis, that is the discontinuity, makes it 
difficult to find a structured and widely applicable ways of fitting the models 
to the data. 

With appropriate measurements, is easy to decide whether the observed bounce 
trajectory entered roll at some point or not. One of the most challenging tasks 
of the modelling will thus be to solve that problem in the initial condition 
space. That is, given the data, future work will be aimed at identifying the 
lower dimensional divide of the space which will identify the initial condition 
that will lead to the bounce with slip only and those that will enter gripping 
at some point.

The aim of identifying such key features is to allow for more precise 
modelling of the compliant surface, in particular in the game of golf. A 
particular basis of nonlinear functions can be selected and thus fitted with 
appropriate parameters using the data available. Avoiding high-dimensional 
problems or finite element solutions remains an open problem. 

A part of the problem at hand is identifying parameters together with their 
physical meaning and the ways to measure them. The golf industry currently 
tends to quantify the ground's firmness and stiffness using either USGA TruFirm 
(Turf-thumper) device or the Clegg Hammer -- the measurements performed by 
these are centred around measuring the deceleration of a free-falling object. 
There is some evidence, however, showing that behaviour of a ball on two turfs 
with similar measurements can be significantly different, thus suggesting that 
important properties are not captured by either of the tools \cite{stri}.

Some immediate extensions of the problem considered are evident and present an 
interesting challenge. One could possibly resolve more complex contacts such as
Hertzian using the current formulation therefore extending the current study to 
the models of elastic half spaces \cite{Barber,haake_apparatus}. A further idea 
to explore would also involve a side spin, that is a spin o . This could 
present a generalisation of the problem to many other disciplines, applicable 
to a further ball-sports but also a wider set of impact-problems in general. 
Both of the aforementioned extensions would require the extension of modelling 
friction in 3D, a problem that presents a far greater order of complexity and 
requires a more considerate approach, whether through the extension of the 
current models or the regularisation of such in higher dimensions -- see 
e.g.~\cite{Antali,cheesman}.


\section*{Acknowledgements}
This work is being partially funded by the EPSRC and by the 
industrial partner R~\&~A Rules 
Ltd., who have also actively participated and facilitated this research. We 
thank the industrial partner for allowing us to run the experimental campaign 
using their research facilities, and particularly acknowledge helpful advice 
from Kristian Jones, Andrew Johnson and Steve Otto. 

\noindent Further thanks go to Mike Jeffrey (University of Bristol) for useful 
discussions on non-smooth  dynamics.

\bibliographystyle{plain}
\bibliography{refs}

\end{document}